\documentclass[a4paper,12pt]{article}
\usepackage{amsmath}
\usepackage{latexsym}
\usepackage{color}
\numberwithin{equation}{section}

\def\R{{\bf R}}

\def\N{{\bf N}}

\def\d{\displaystyle}
\def\e{{\varepsilon}}

\def\wt{\widetilde}

\def\p{\partial}

\newtheorem{thm}{Theorem}[section]

\newtheorem{lem}{Lemma}[section]
\newtheorem{prop}{Proposition}[section]

\title{The lifespan of solutions of semilinear
wave equations with the scale-invariant damping in two space
dimensions}
\author{
Takuto Imai
\footnote{Accenture Japan Ltd,
Harumi Triton Square Office Tower Z,
1-8-12 Harumi, Chuo-ku, Tokyo, 104-0053, Japan.
e-mail: takuto.imai@accenture.com.}
, Masakazu Kato
\footnote{College of Liberal Arts, Mathematical Science Research Unit, 
Muroran Institute of Technology, 27-1 Mizumoto-cho, 
Muroran, Hokkaido 050-8585, Japan.
email: mkato@mmm.muroran-it.ac.jp.}
, Hiroyuki Takamura
\footnote{Mathematical Institute,
Tohoku University,
 Aoba, Sendai 980-8578, Japan.
e-mail: hiroyuki.takamura.a1@tohoku.ac.jp.}\\
and
Kyouhei Wakasa
\footnote{
Department of Creative Engineering, National Institute of Technology, Kushiro College, 2-32-1 Otanoshike-Nishi, Kushiro-Shi, Hokkaido 084-0916, Japan.
e-mail: wakasa@kushiro-ct.ac.jp.
}}
\date{
\[
\begin{array}{ll}
\mbox{\footnotesize{\bf Keywords:}}
& \mbox{\footnotesize semilinear wave equation, scale-invariant damping, lifespan}\\
\mbox{\footnotesize{\bf MSC2010:}}
& \mbox{\footnotesize primary 35L71, secondary 35B44}\\
\end{array}
\]
}
\begin{document}
\maketitle
\begin{abstract}
In this paper, we study the initial value problem for semilinear wave equations 
with the time-dependent and scale-invariant damping in two dimensions. 
Similarly to the one dimensional case by Kato, Takamura and Wakasa in 2019,
we obtain the lifespan estimates of the solution for a special constant in the damping term,
which are classified by total integral of the sum of the initial position and speed.
\par
The key fact is that, only in two space dimensions,
such a special constant in the damping term is a threshold
between \lq\lq wave-like" domain and \lq\lq heat-like" domain.
As a result, we obtain a new type of estimate especially for the critical exponent.
\end{abstract}

\section{Introduction}
\par
We are concerned with the following initial value problem
for semilinear wave equations with the scale-invariant damping:
\begin{equation}
\label{IVP1}
\left\{
\begin{array}{ll}
\d v_{tt}-\Delta v+\frac{\mu}{1+t}v_t=|v|^p
&\mbox{in}\quad \R^n\times[0,\infty),\\
v(x,0)=\e f(x),\ v_t(x,0)=\e g(x),
& x\in\R^n,
\end{array}
\right.
\end{equation}
where $v=v(x,t)$ is a real valued unknown function,
$\mu>0$, $p>1$, $n\in\N$, the initial data $(f,g) \in H^1(\R^n) \times L^2(\R^n)$ has compact support,
and $\e>0$ is a \lq\lq small" parameter.
\par
It is interesting to look for the critical exponent $p_c(n)$ such that
\begin{align*}
\left\{
\begin{array}{ll}
	p> p_{c}(n) \ (\mbox{and may have an upper bound}) & 
	\Longrightarrow \ T(\e)=\infty,\\
	1 < p \leq p_{c}(n) & \Longrightarrow \ T(\e) < \infty,
\end{array}
\right.
\end{align*}
where $T(\e)$ is the lifespan, the maximal existence time,
of the energy solution of (\ref{IVP1}) with an arbitrary fixed non-zero data.
Then, we have the following conjecture:
\begin{align}
\left\{
\begin{array}{lll}
	\mu > \mu_{0}(n) & \Longrightarrow \ p_{c}(n)=p_{F}(n) 
	&\quad (\mbox{heat-like}),\\
	\mu = \mu_{0}(n) & \Longrightarrow \ p_{c}(n)=p_{F}(n)=p_{S}(n+\mu) 
	&\quad (\mbox{intermediate}),\\
	0< \mu < \mu_{0}(n) & \Longrightarrow \ p_{c}(n)=p_{S}(n+\mu)
	&\quad (\mbox{wave-like}),
\end{array}
\right.\label{dai1.1}
\end{align}
where
\begin{align}
\label{mu_0}
	\mu_{0}(n):=\frac{n^2+n+2}{n+2}.
\end{align}
Here
\begin{align}
\label{fujita}
	p_{F}(n):=1+\frac{2}{n}
\end{align}
is the so-called Fujita exponent which is the critical exponent of the associated semilinear
heat equations $v_{t}-\Delta v=v^p$,
and
\begin{align}
\label{strauss}
	p_{S}(n):=
	\left\{
	\begin{array}{ll}
		\infty & (n=1),\\
		\d\frac{n+1+\sqrt{n^2+10n-7}}{2(n-1)} & (n \geq 2)
	\end{array}
	\right.
\end{align}
is the so-called  Strauss exponent which is the critical exponent of the associated
semilinear wave equations $v_{tt}-\Delta v= |v|^p$.
We note that $p_{S}(n)$ $(n \geq 2)$ is the positive root of
\begin{align}
\label{gamma}
	\gamma(p,n):=2+(n+1)p-(n-1)p^2=0
\end{align}
and $0 < \mu < \mu_0(n)$ is equivalent to $p_{F}(n) < p_{S}(n+\mu)$.
\par
The conjecture (\ref{dai1.1}) shows the critical situation of our problem in the following sense.
If one replaces the damping term $\mu v_t/(1+t)$ in (\ref{IVP1}) by $\mu v_t/(1+t)^\beta$,
then one can see that there is no such a $p_c(n)$, namely $T(\e)=\infty$ for any $p>1$ when $\beta<-1$,
the so-called over damping case.
Moreover one has $p_c(n)=p_F(n)$ for any $\mu>0$ when $-1\le\beta<1$,
the so-called effective damping case, and
$p_c(n)=p_S(n)$ for any $\mu>0$ (it can be any $\mu\in\R$) when $\beta>1$,
the so-called scattering damping case.
Therefore one may say that the so-called scale-invariant case, $\beta=1$,
is an intermediate situation between wave-like,
in which the critical exponent is related to $p_S(n)$,
and heat-like, in which the critical exponent is $p_F(n)$.
To see all the references above results, for example,
see introductions of related papers to the scattering damping case,
Lai and Takamura \cite{LT18} (the sub-critical case),
Wakasa and Yordanov \cite{WY19} (the critical case),
Liu and Wang \cite{LW} (partial result of the super-critical case).
\par
For the conjecture (\ref{dai1.1}), D'Abbicco \cite{D15} has obtained the heat-like existence
partially with
\begin{align*}
\mu \geq \left\{
\begin{array}{ll}
	5/3 & \mbox{for} \ n=1\ \mbox{(cf. $\mu_0(1)=4/3$)},\\
	3 & \mbox{for} \ n=2\ \mbox{(cf. $\mu_0(2)=2$)},\\
	n+2 & \mbox{for} \ n \geq 3,
\end{array}
\right.
\end{align*}
while Wakasugi \cite{Wakasugi14} has obtained the blow-up parts
in $1<p<p_{F}(n)$ for $\mu \geq 1$
and $1<p<p_{F}(n+\mu-1)$ for $0< \mu < 1$.
In the damped case $\mu >0$, his second result is the first blow-up result for super-Fujita exponents
which are larger than $p_{F}(n)$.

\par
In this paper, we consider a special case of $\mu=2$. 
The speciality of this value is clarified by setting
\[
u(x,t):=(1+t)^{\mu/2}v(x,t),
\]
where $v$ is the solution to (\ref{IVP1}).
Then, $u$ satisfies
\begin{equation}\label{IVP2}
\left\{
\begin{array}{l}
\d u_{tt}-\Delta u+\frac{\mu(2-\mu)}{4(1+t)^2}u=\frac{|u|^p}{(1+t)^{\mu(p-1)/2}}
\quad\mbox{in}\quad \R^n\times[0,\infty),\\
u(x,0)=\e f(x),\ u_t(x,0)=\e \{\mu f(x)/2+g(x)\},\quad x\in\R^n,
\end{array}
\right.
\end{equation}
so that all the technics in the analysis on semilinear wave equations can be employed
and we may discussed about not only the energy solution but also the classical solution.
In fact, via this reduced problem (\ref{IVP2}),
D'Abbicco, Lucente and Reissig \cite{DLR15} have proved
the intermediate part of the conjecture (\ref{dai1.1}) for $n=2$ and the wave-like part for $n=3$ when $\mu=2$.
We note that the assumption of the radial symmetry is considered in \cite{DLR15}
for the existence part in $n=3$.
Moreover, D'Abbicco and Lucente \cite{DL15} have obtained the wave-like existence part of (\ref{dai1.1}) for odd $n \geq 5$ when $\mu=2$ also with radial symmetry.
 \par
 In the case of $\mu \neq 2$, Lai, Takamura and Wakasa 
 \cite{LTW17} have 
 first studied the wave-like blow-up of the conjecture (\ref{dai1.1}) with a loss replacing
 $\mu$ by $\mu/2$ in the sub-critical case.
 Initiating this result, Ikeda and Sobajima \cite{IkeSoba} have obtained the blow-up part of (\ref{dai1.1}).
\par
For the lifespan estimate, one may expect that
\begin{align}
 T(\e) \sim \left\{
 \begin{array}{ll}
 C \e^{-(p-1) / \left\{ 2-n(p-1)\right\}} & \mbox{for} \ 1<p<p_{F}(n),\\
 \exp \left( C \e^{-(p-1)}\right) & \mbox{for} \ p=p_{F}(n)
 \end{array}
 \right.\label{dai1.2}
\end{align}
for the heat-like domain $\mu > \mu_{0}(n)$ and 
\begin{align}
 T(\e) \sim \left\{
 \begin{array}{ll}
 C \e^{-2p(p-1) / \gamma(p,n+\mu)} & \mbox{for} \ 1<p<p_{S}(n+\mu),\\
 \exp \left( C \e^{-p(p-1)}\right) & \mbox{for} \ p=p_{S}(n+\mu)
 \end{array}
 \right.\label{dai1.3}
\end{align}
for the wave-like domain $0< \mu < \mu_{0}(n)$.
Recall the definitions of $\mu_0(n)$, $p_F(n)$, $p_S(n)$ and $\gamma(p,n)$
in (\ref{mu_0}), (\ref{fujita}), (\ref{strauss}) and (\ref{gamma}).
Here $T(\e) \sim A(\e,C)$
stands for the fact that there are positive constants, $C_{1}$ and $C_{2}$,
independent of $\e$ satisfying 
$A(\e,C_1) \leq T(\e) \leq A(\e,C_2)$.
Actually, (\ref{dai1.2}) for $n=1$ and
$\mu=2>\mu_0(1)=4/3$ is obtained by Wakasa \cite{Wa16}, and (\ref{dai1.3}) is obtained by Kato and Sakuraba \cite{KS} for $n=3$ and $\mu=2<\mu_0(3)=14/5$.
One may refer Lai \cite{Lai} for the existence part of weaker solution.
Moreover, the upper bound of (\ref{dai1.2}) in the sub-critical case is obtained by Wakasugi \cite{Wakasugi14}. Also the upper bound of (\ref{dai1.3}) is 
obtained by Ikeda and Sobajima \cite{IkeSoba} in the critical case, later it is reproved
 by Tu and Li \cite{TL2}, and Tu and Li \cite{TL1} in the sub-critical case.
\par
In the non-damped case of $\mu=0$, it is known that (\ref{dai1.3}) is true for $n \geq 3$, or $p>2$ and $n=2$, The open part around this fact is $p=p_{S}(n)$ for $n \geq 9$.
In other cases, (\ref{dai1.3}) is still true if $\int_{\R^n} g(x) dx=0$. On the other hand,
we have
\begin{align}
	T(\e) \sim \left\{
	\begin{array}{ll}
	C \e^{-(p-1)/2} & \mbox{for} \ n=1,\\
	C \e^{-(p-1)/(3-p)} & \mbox{for} \  n=2 \ \mbox{and} \ 1<p<2,\\
	C a(\e) & \mbox{for} \ n=2 \ \mbox{and} \ p=2
	\end{array}
	\right.\label{dai1.4}
\end{align}
if $\int_{\R^n} g(x) dx \neq 0$, where $a=a(\e)$ is a positive number satisfying
$\e^2 a ^2 \log(1+a)=1$. We note that the bounds in (\ref{dai1.4}) are smaller than
the one of the first line in (\ref{dai1.3}) with $\mu=0$ in each case. For all the references in the case of $\mu=0$, 
see Introduction of Imai, Kato, Takamura and Wakasa \cite{IKTW}.
\par
The remarkable fact is that even if $\mu$ is in the heat-like domain, the lifespan estimate for (\ref{IVP1}) is similar to the one for non-damped case.
Indeed, for $n=1$ and $\mu=2>\mu_{0}(1)=4/3$, Kato, Takamura and Wakasa \cite{KTW} 
show that the result on (\ref{dai1.2}) by Wakasa \cite{Wa16} mentioned above is true only
if $\int_{\R} \{ f(x)+g(x) \} dx \neq 0$. More precisely, they have obtained that
\begin{align}
	T(\e) \sim \left\{
	\begin{array}{ll}
		C \e^{-2p(p-1)/\gamma(p,3)} & \mbox{for} \ 1<p<2,\\
		C b(\e) & \mbox{for} \ p=2,\\
		C \e^{-p(p-1)/(3-p)}  & \mbox{for} \ 2<p<3,\\
		\exp(C \e^{-p(p-1)}) & \mbox{for} \ p=p_{F}(1)=3,
	\end{array}
	\right.\label{dai1.5}
\end{align}
if $\int_{\R} \{ f(x)+g(x) \} dx=0$, where $b=b(\e)$ is a positive number
satisfying $\e^2 b \log(1+b)=1$. We note that the bounds in (\ref{dai1.5}) are
 larger than those in (\ref{dai1.2}) with $n=1$
and $\mu=2$ in each case.

\par
Our aim in this paper is to show the lifespan estimates for (\ref{IVP1})
in two dimensional case, $n=2$, with $\mu=2$
which is similar to one dimensional case as above.
We note $p_c(2)=p_F(2)=p_S(2+2)=2$ and $\mu_0(2)=2$.
More precisely, we shall show that
\begin{equation}
\label{dai1.6}
	T(\e)\sim\left\{
	\begin{array}{ll}
	c\e^{-(p-1)/(4-2p)}
	&\mbox{for} \quad\d 1<p<2,\\
	\exp(c\e^{-1/2})
	& \mbox{for}\quad\d p=2
\end{array}
\right.
\end{equation}
if $\int_{\R^2} \{ f(x)+g(x)\} dx \neq 0$,
and
\begin{equation}
\label{dai1.7}
T(\e)\sim\left\{
	\begin{array}{ll}
	c\e^{-2p(p-1)/\gamma(p,4)}
	&\mbox{for} \quad\d 1<p<2,\\
	\exp(c\e^{-2/3})
	& \mbox{for}\quad\d p=2
\end{array}
\right.
\end{equation}
if $\int_{\R^2} \left\{ f(x)+g(x) \right\} dx = 0$.
We note that the critical cases in (\ref{dai1.6}) and (\ref{dai1.7}) are new
in the sense that they are different from (\ref{dai1.2}) and (\ref{dai1.3}).
(\ref{dai1.6}) and (\ref{dai1.7}) are announced in Introduction of Lai and Takamura \cite{LT18},
but there are typos in the exponents of $\e$ in the critical case.
\par
The strategy of proofs in this paper is based on point-wise estimates of the solution.
In the existence part,
we employ the classical iteration argument for semilinear wave equations without damping term,
which is first introduced by John \cite{John79} in three space dimensions,
and its variant, which is developed by Imai, Kato, Takamura and Wakasa \cite{IKTW}
in two space dimensions.
In the blow-up part, we also employ an improved version of Kato's lemma
on ordinary differential inequality by Takamura \cite{Takamura15}
for the sub-critical cases.
We note that, till now, the so-called test function method
such as in Ikeda, Sobajima and Wakasa \cite{ISW} cannot be applicable to delicate analysis
to catch the logarithmic growth of the solution in the case of $p=2$ in (\ref{dai1.4}), (\ref{dai1.5}), (\ref{dai1.6}) and (\ref{dai1.7}).
Therefore we employ the so-called slicing method of the blow-up domain for the critical case,
which is introduced by Agemi, Kurokawa and Takamura \cite{AKT00}
to handle weakly coupled systems of non-damped semilinear wave equations with critical exponents.
\par
This paper is organized as follows.
In the next section, our goals, (\ref{dai1.6}) and (\ref{dai1.7}), are described in four theorems,
and we introduce the linear decay estimate and basic lemmas for a-priori estimates.
Section 3, or Section 4, is devoted to the proof of the lower bound, or upper bound,
of the lifespan respectively.


\section{Theorems and preliminaries}
\label{section:Pre}
\par
In this section, we state our results (\ref{dai1.6}) and (\ref{dai1.7}) in four theorems.
After them, we list useful point-wise estimates of linear wave equations.
For the sake of the simplicity, we assume that
\begin{align}
	\mbox{supp} (f,g) \subset 
	\{ x \in \R^2 \ : \ |x| \leq k \},\ k\ge1
\label{supp_u}
\end{align}
throughout this paper.
\par
The existence parts of our goals in (\ref{dai1.6}) and (\ref{dai1.7})
are guaranteed by the following two theorems.
Recall the definitions of $\mu_0(n)$, $p_F(n)$, $p_S(n)$ and $\gamma(p,n)$
respectively in (\ref{mu_0}), (\ref{fujita}), (\ref{strauss}) and (\ref{gamma}).

\begin{thm}
\label{T.1.1}
Let $n=2$, $\mu=\mu_0(2)=2$ and $1<p\le p_F(2)=p_{S}(4)=2$.
Assume that $(f,g)\in C_0^3(\R^2)\times C_0^2(\R^2)$ satisfies (\ref{supp_u}).
Then, there exists a positive constant $\e_0=\e_0(f,g,p,k)$ such that 
(\ref{IVP2}) admits a unique solution $u\in C^2(\R^2\times[0,T))$
if $p=2$, or the integral equation associated with (\ref{IVP2}) admits a unique solution $u\in C^1(\R^2\times[0,T))$ otherwise, 
as far as $T$ satisfies
\[
	T\le\left\{
	\begin{array}{ll}
	c\e^{-(p-1)/(4-2p)}
	&\mbox{if} \quad\d 1<p<2,\\
	\exp(c\e^{-1/2})
	& \mbox{if}\quad\d p=2
\end{array}
\right.
\]
for $0<\e\le\e_0$, where $c$ is a positive constant independent of $\e$. 
\end{thm}

\begin{thm}
\label{T.1.2}
Suppose that the assumptions in Theorem \ref{T.1.1} are fulfilled.
Assume additionally that 
\[
\int_{\R^2}\{f(x)+g(x)\}dx=0.
\]
Then, there exists a positive constant $\e_0=\e_0(f,g,p,k)$ such that 
(\ref{IVP2}) admits a unique solution $u\in C^2(\R^2\times[0,T))$
if $p=2$, or the integral equation associated with (\ref{IVP2}) admits a unique solution $u\in C^1(\R^2\times[0,T))$ otherwise, 
as far as $T$ satisfies
\[
	T\le\left\{
	\begin{array}{ll}
	c\e^{-2p(p-1)/\gamma(p,4)}
	&\mbox{if} \quad\d 1<p<2,\\
	\exp(c\e^{-2/3})
	& \mbox{if}\quad\d p=2
\end{array}
\right.
\]
for $0<\e\le\e_0$, where $c$ is a positive constant independent of $\e$.
\end{thm}

On the other hand, the blow-up parts of our goals in (\ref{dai1.6}) and (\ref{dai1.7})
are guaranteed by the following two theorems.

\begin{thm}
\label{T.2.1}
Let $n=2$, $\mu=\mu_0(2)=2$, $1<p\le 2=p_F(2)=p_S(4)$.
Assume that $(f,g) \in  C_0^2(\R^2) \times C_0^1(\R^2)$ satisfies (\ref{supp_u}).
Suppose that the integral equation associated with (\ref{IVP2}) has a solution 
$u\in C^1(\R^2 \times[0,T))$ with {\rm supp} $u\subset \{(x,t) \in \R^2 \times
[0,\infty) : |x| \leq t+k \}$. 
Then, there exists a positive constant $\e_1=\e_1(f,g,p,k)$ such that  the solution cannot exist whenever $T$ satisfies
\begin{equation*}
T\ge\left\{
	\begin{array}{ll}
	c\e^{-(p-1)/(4-2p)}
	&\mbox{if} \quad\d 1<p<2,\ f(x) \equiv 0 \ \mbox{and} \
	g(x) \geq 0 \ (\not\equiv0),  \\
	\exp(c\e^{-1/2})
	& \mbox{if}\quad\d p=2 \  \mbox{and} \ \int_{\R^2} \{ f(x) +g(x) \} dx >0
\end{array}
\right.
\end{equation*}
for $0<\e\le\e_1$, where $c$ is a positive constant independent of $\e$.
\end{thm}

\begin{thm}
\label{T.2.2}
Suppose that the assumptions in Theorem \ref{T.2.1} are fulfilled.
Assume additionally that
\[
f(x)+g(x)\equiv0.
\]
Then, there exists a positive constant $\e_1=\e_1(f,g,p,k)$ such that the solution 
cannot exist whenever $T$ satisfies
\begin{equation*}
T\ge\left\{
	\begin{array}{ll}
	c\e^{-2p(p-1)/\gamma(p,4)}
	&\mbox{if} \quad\d 1<p<2 \ \mbox{and} \ f(x)\ge0 \ (\not\equiv0),\\
	\exp(c\e^{-2/3})
	& \mbox{if}\quad\d p=2 \ \mbox{and} \ \int_{\R^2} f(x) dx <0
\end{array}
\right.
\end{equation*}
for $0<\e\le\e_1$, where $c$ is a positive constant independent of $\e$.
\end{thm}

\par
From now on, we introduce some definitions and useful lemmas.
For $(x,t) \in \R^2 \times [0,\infty)$, we set
 \begin{equation}
\label{u_L}
\begin{array}{ll}
\d u_L(x,t):=\frac{\partial}{\partial t}R(f|x,t)+R(f+g|x,t),\\
\d R(\phi|x,t):=\frac{1}{2\pi}\int_{|x-y|\le t}\frac{\phi(y)}{\sqrt{t^2-|x-y|^2}}dy
=\frac{t}{2\pi}\int_{|\xi|\le1}\frac{\phi(x+t\xi)}{\sqrt{1-|\xi|^2}}d\xi.
\end{array}
\end{equation}
When $(f,g) \in C_{0}^3(\R^2) \times C_{0}^2(\R^2)$, we note that $u_L$ satisfies that
\[
\left\{
\begin{array}{ll}
\d (u_L)_{tt}-\Delta u_L=0
& \mbox{in}\ \R^2\times[0,\infty),\\
u_L(x,0)=f(x),\ (u_L)_t(x,0)=f(x)+g(x),
& x\in \R^2
\end{array}
\right.
\]
in the classical sense, and it holds
\[
{\rm supp}\ u_L\subset\{(x,t)\in\R^2\times[0,\infty)\ :\ |x|\le t+k\}.
\]
We introduce the decay estimates for the solutions of (\ref{u_L}) which will be used 
in the proof of Theorem \ref{T.1.1} and Theorem \ref{T.1.2}.
For the proof, see Lemma 2.1 in \cite{IKTW}.

\begin{lem}[Imai, Kato, Takamura and Wakasa \cite{IKTW}]
\label{lem:decay_est_v}
Let $u_L$ be the one in (\ref{u_L}).
Then, there exist positive constants
$\wt{C_0}=\wt{C_0}(k)$ and
$C_{0}=C_0(\|f\|_{W^{3,1}(\R^2)},\|g\|_{W^{2,1}(\R^2)},k)$ 
such that $u_L$ satisfies
\[
\begin{array}{lllll}
&\d \sum_{|\alpha|\le1}|\nabla_x^\alpha u_L(x,t)|\\
&\le\d \left|\int_{\R^2}\{f(x)+g(x)\}dx\right|\cdot\frac{\wt{C_{0}}}{(t+|x|+2k)^{1/2}(t-|x|+2k)^{1/2}}\\
&\d \quad+\frac{C_0}{(t+|x|+2k)^{1/2}(t-|x|+2k)^{3/2}}
\end{array}
\]
in $\R^2\times[0,\infty)$.
\end{lem}

\par
Next, we prepare the following decay estimate which 
will be employed in the proof of Theorem \ref{T.2.1} and Theorem \ref{T.2.2}.

\begin{lem}
\label{lem4.5}
Let $u_{L}$ be the one in (\ref{u_L}).
For $t-|x| \geq 2k$ and $t \geq 4k$, there exists  a positive constant
$C=C(f,g)$ such that
\begin{align}
	\left| u_{L}(x,t)-\frac{\d \int_{\R^2} \{ f(x)+g(x) \}  dx }{2 \pi (t+|x|)^{1/2}(t-|x|)^{1/2}} \right|
	\leq\frac{C k}{(t+|x|)^{1/2} (t-|x|)^{3/2}},
	\label{dai4.50}
\end{align}
moreover that
\begin{align}
\left| u_{L}(x,t)+\frac{ t \d\int_{\R^2} f(x)  dx }{2 \pi (t+|x|)^{3/2}(t-|x|)^{3/2}} \right|
\leq\frac{C k}{(t+|x|)^{1/2} (t-|x|)^{5/2}}
\label{dai4.51}
\end{align}
provided $f(x)+g(x)\equiv0$.
\end{lem}

\par\noindent
{\bf Proof.}
First we shall prove (\ref{dai4.51}).
Denote $r:=|x|$.
For $t-r \ge 2k$ and $t \ge 4k$, we shall split the domain into the interior domain $t \ge 2r$ and the exterior domain $t \le 2r$. We set
\begin{align}
	&D_{int}:=\{(x,t) \in \R^2 \times [0,\infty) \ : \ t \geq 2r,\ t \geq 4k \},
	\nonumber\\
	&D_{ext}:=\left\{ (x,t) \in \R^2 \times [0,\infty) \ : \ 
	r+2k \leq t \leq 2r \right\}.\nonumber
\end{align}
\par
First, we prove (\ref{dai4.51}) in $D_{int}$.
Since $f(x)+g(x)\equiv0$, $(\ref{supp_u})$ and 
\begin{align}
 	|x-y|\leq r+|y| \leq \frac{t}{2}+k \leq t \quad \mbox{for} \ (x,t) \in D_{int} \
 	\mbox{and} \ |y| \leq k, \nonumber
\end{align}
we can rewrite $u_{L}(x,t)$ in (\ref{u_L}) as
\[
	u_{L}(x,t) = \frac{\partial}{\partial t}
	\left\{ 
	\frac{1}{2 \pi} \int_{\R^2} \frac{f(y)}{\sqrt{t^2-|x-y|^2}} dy
	\right\}\nonumber\\
	=-\frac{t}{2 \pi} \int_{\R^2} \frac{f(y)}{(t^2-|x-y|^2)^{3/2}} dy.\nonumber
\]
This expression gives us
\begin{align}
\label{dai4.52}
	&\left|  
	\frac{2\pi}{t} u_{L}(x,t) + \frac{1}{(t^2-r^2)^{3/2}} \int_{\R^2} f(y) dy
	\right|\nonumber\\
	&\leq \frac{1}{(t^2-r^2)^{3/2}} \int_{\R^2} \frac{|h_{1}(x,y,t)|}{(t^2-|x-y|^2)^{3/2}}
	|f(y) | dy,
\end{align}
where 
\begin{align}
	h_{1}(x,y,t):=(t^2-r^2)^{3/2}-(t^2-|x-y|^2)^{3/2}.\nonumber
\end{align}
Using the Taylor expansion with respect to $y$ at the origin, we get
\begin{align}
	h_{1}(x,y,t)=3(t^2-|x-\theta y|^2)^{1/2} 
	\left\{ -<x,y> + \theta |y|^2 \right\} 
		\label{dai2.10}
\end{align}
with $0<\theta<1$.
For $(x,t) \in D_{int}$ and $|y| \le k$, we obtain
\begin{align}
	(t^2-|x-\theta y|^2)^{1/2} &\leq t+ r+|y|
	\leq \frac{3}{2}t+k\label{dai2.11}
\end{align}
and
\begin{align}
	|-<x,y>+ \theta |y|^2| \leq\left(\frac{t}{2}+k\right)k.\label{dai2.12}
\end{align}
From (\ref{dai2.10}), (\ref{dai2.11}) and (\ref{dai2.12}), it follows that
\begin{align}
	|h_{1}(x,y,t)| \leq 3 \left(\frac{3}{2}t+k\right)
	\left(\frac{t}{2}+k\right)k\leq C k t^2.\label{dai4.53}
\end{align}
Therefore, combining (\ref{dai4.52}), (\ref{dai4.53}) and
\begin{align}
	t+|x-y| \geq t-|x-y| \geq t-r-\frac{t}{4}
	\geq \frac{t-r}{2},\nonumber
\end{align}
we have 
\begin{align}
	\left| \frac{2 \pi}{t}u_{L} +\frac{1}{(t^2-r^2)^{3/2}} \int_{\R^2} f(y) dy\right|
	&\leq \frac{C k t^2}{(t+r)^{3/2}(t-r)^{9/2}}.\nonumber
\end{align}
Since $3(t-r) \geq t+r$ holds in $D_{int}$, we obtain (\ref{dai4.51}) in $D_{int}$.

\par
Next, we prove (\ref{dai4.51}) in $D_{ext}$.
Here, we employ the following different representation formula from (\ref{u_L})
which is established by (6.24) in H\"ormander \cite{H97}:
\begin{align}
	u_{L}(x,t)=\frac{\partial}{\partial t}
	\left(
	\frac{1}{2\sqrt{2} \pi \sqrt{r}} \int_{\rho-\rho^2z/2}^{\infty}
	\frac{I(f)(s,\omega,z)}{\sqrt{s-\rho+\rho^2z/2}}ds
	\right),\label{dai4.54}
\end{align}
where $\omega=x/r \in S^{1}$, $\rho=r-t$, $z=1/r$ and
\[
I(f)(s,\omega,z):=
\int_{s=<\omega,y>-|y|^2z/2} f(y) dS_{y}.
\]
For $(x,t) \in D_{ext}$ and $|y| \leq k$, we have
\begin{align}
	\left|<\omega,y> -\frac{|y|^2z}{2}\right| &\leq |y| + \frac{|y|^2z}{2}
	\leq \frac{5k}{4}.\label{dai4.55}
\end{align}
Since
\[
\rho-\frac{\rho^2 z}{2}=-\frac{(t+r)(t-r)}{2r} \leq -2k,
\ -1+\rho z=-\frac{t}{r}
\ \mbox{and}\ \frac{\partial}{\partial t}=-\frac{\partial}{\partial \rho},
\]
it follows from (\ref{dai4.54}) and (\ref{dai4.55}) that
\begin{align}
	u_{L}&=\frac{\partial}{\partial t}
	\left(
	\frac{1}{2 \sqrt{2} \pi \sqrt{r}}
	\int_{-5k/4}^{5k/4}
	\frac{I(f)(s,\omega,z)}{\sqrt{s-\rho+\rho^2z/2}}ds
	\right)\nonumber\\
	&=-\frac{t}{4 \sqrt{2} \pi r^{3/2}}
	\int_{-5k/4}^{5k/4}
	\frac{I(f)(s,\omega,z)}{(s-\rho+\rho^2z/2)^{3/2}}ds.\nonumber
\end{align}
Hence, we obtain
\begin{align}
	&\left| \frac{4\sqrt{2} \pi r^{3/2}}{t} u_{L}+
	\left\{ \frac{2r}{(t+r)(t-r)} \right\}^{3/2}
	\int_{\R^2} f(y) dy \right|\nonumber\\
	&\leq
	\left\{ \frac{2r}{(t+r)(t-r)} \right\}^{3/2}
	\int_{-5k/4}^{5k/4}
	\frac{|h_2(\rho,s,z)|I(|f|)(s,\omega,z)}{(s-\rho+\rho^2z/2)^{3/2}}ds,\label{dai4.56}
\end{align}
where
\begin{align}
	h_{2}(\rho, s, z):=\left(-\rho+\frac{\rho^2z}{2}\right)^{3/2}
	-\left(s-\rho+\frac{\rho^2z}{2}\right)^{3/2}.\label{dai4.57}
\end{align}
Making use of the Taylor expansion with respect to $s$ at the origin, we have
from (\ref{dai4.57})
\begin{align}
	h_{2}(\rho,s,z)=-\frac{3}{2} \left(\theta s -\rho+\frac{\rho^2z}{2}\right)^{1/2} s
	\quad \mbox{with} \ 0 < \theta < 1.\label{dai2.13}
\end{align}
Since $\rho=r-t$ and $z=1/r$, for $|s| \le 5k/4$, we obtain
\begin{align}
	\left|\theta s -\rho + \frac{\rho^2 z}{2}\right|^{1/2}|s|
	&\leq \left\{ \frac{5k}{4}+ \frac{(t+r)(t-r)}{2r} \right\}^{1/2}
	\cdot \frac{5}{4}k\nonumber\\
	&\leq \frac{C k (t+r)^{1/2}(t-r)^{1/2}}{r^{1/2}}\label{dai4.58}
\end{align}
and
\begin{align}
	s-\rho +\frac{\rho^2z}{2} 
	\geq -\frac{5k}{4}+t-r
	\geq \frac{3}{8}(t-r).\label{dai4.59}
\end{align}
Hence, for $(x,t) \in D_{ext}$, it follows from (\ref{dai4.56}), (\ref{dai2.13}),
(\ref{dai4.58}) and (\ref{dai4.59}) that
\begin{align}
	&\left| \frac{4 \sqrt{2} \pi r^{3/2}}{t} u_{L} + 
	\left\{ \frac{2r}{(t+r)(t-r)} \right\}^{3/2}
	\int_{\R^2} f(y) dy \right|\nonumber\\
	&\leq \frac{C k r}{(t+r)(t-r)^{5/2}} \int_{\R^2} |f(y)| dy.\nonumber
\end{align}
Therefore, we obtain (\ref{dai4.51}) in $D_{ext}$.
\par
Finally, we show (\ref{dai4.50}).
It follows from the proof of Lemma 2.1 in \cite{IKTW} that 
\begin{align}
	\left| R(f+g|x,t) -\frac{\d\int_{\R^2} \left\{  f(x)+g(x) \right\} dx}{2 \pi (t+r)^{1/2}(t-r)^{1/2}} \right|
	\leq \frac{C k}{(t+r)^{1/2} (t-r)^{3/2}}\label{dai4.60}
\end{align}
for $t-r \ge 2k$ and $t \ge 4k$, where $R(f+g|x,t)$ is defined in (\ref{u_L}).
From (\ref{u_L}), (\ref{dai4.51}) and $t-r \ge 2k$, we have
\begin{align}
	 \left| \frac{\partial}{\partial t} R(f | x,t) \right| 
	\le \frac{Ck}{(t+r)^{1/2}(t-r)^{3/2}}.\label{dai4.60.1}
\end{align}
From (\ref{u_L}), (\ref{dai4.60}) and (\ref{dai4.60.1}), we obtain (\ref{dai4.50}).
This completes the proof.
\hfill$\Box$
\vskip10pt

\par
In what follows, we consider the following integral equation:
\begin{align}
	u(x,t)=u^{0}(x,t)+L(F)(x,t) \quad \mbox{for} \ (x,t) \in \R^2 \times [0,\infty),
	\label{dai2.7}
\end{align}
where we set
\begin{align}
	u^{0}:= \e u_L\label{dai2.7.1}
\end{align}
and 
\begin{align}\label{L}
	L(F)(x,t):=\frac{1}{2 \pi} \int_{0}^{t}\frac{t-\tau}{(1+\tau)^{p-1}}d \tau 
	\int_{|\xi| \leq 1} \frac{F(x+(t-\tau)\xi,\tau)}{ \sqrt{1 - |\xi|^2}}
	 d \xi 
\end{align}
for $F\in C(\R^2\times[0,\infty))$.
We note that (\ref{L}) solves
\[
\left\{
\begin{array}{ll}
u_{tt}-\Delta u=(1+t)^{-(p-1)}F & \mbox{in}\ \R^2\times[0,\infty),\\
u(x,0)=0,\ u_t(x,0)=0, & x\in\R^2
\end{array}
\right.
\]
when $F\in C^2(\R^2\times[0,\infty))$.
Then, the following lemma is one of the basic tools.

\begin{lem}[Agemi and Takamura \cite{AT92}]\label{lem:int_new_L}
Let $L$ be the linear integral operator defined by {\rm (\ref{L})} 
and $\Psi=\Psi(|x|,t)\in C( \R^2 \times [0,\infty))$.
Then we have
\[
L\left(\Psi\right)(x,t)
=L_{1}\left(\Psi\right)(r,t)+L_{2}\left(\Psi\right)(r,t),\ r=|x|,
\]
where $L_{j}\left(\Psi\right)\ (j=1,2)$ are defined by 
\[
\begin{array}{lll}
L_{1}\left(\Psi\right)(r,t)\\
:=
	\d \frac{2}{\pi}\int_{0}^{t} (1+\tau)^{-(p-1)} d\tau
	\int_{|t-r-\tau|}^{t+r-\tau}
	\lambda \Psi(\lambda,\tau)d\lambda
	\int_{|\lambda-r|}^{t-\tau}
	\frac{\rho h(\lambda,\rho;r)}{\sqrt{(t-\tau)^2-\rho^2}}d\rho,
\end{array}
\]
\[
\begin{array}{llll}
	L_{2}\left(\Psi\right)(r,t)\\
\d:=
	\frac{2}{\pi}\int_{0}^{(t-r)_{+}} (1+\tau)^{-(p-1)} d\tau
	\int_{0}^{t-r-\tau}\!\!\!\lambda \Psi(\lambda,\tau)d\lambda
	\int_{|\lambda-r|}^{\lambda+r}
	\frac{\rho h(\lambda,\rho;r)}{\sqrt{(t-\tau)^2-\rho^2}}d\rho,
\end{array}
\]
where $a_+:=\max\{a,0\}$ and
\begin{align*}
	h(\lambda,\rho;r):=\{ (\rho^2-(\lambda-r)^2)
	((\lambda + r)^2-\rho^2) \}^{-\frac{1}{2}}.
\end{align*}
\par
Moreover, the following estimates hold in $[0,\infty)^2$:
\begin{align}
	|L_{1}\left(\Psi\right)(r,t)| & \le 
	\frac{1}{\sqrt{2}}\int_{0}^{t} (1+\tau)^{-(p-1)} d \tau
	\int_{|r-t+\tau|}^{r+t-\tau}
	\frac{\lambda |\Psi(\lambda,\tau)|d\lambda}
	{(\sqrt{r} \ \mbox{or} \ \sqrt{\lambda})\sqrt{\tau+\lambda-t+r}},\nonumber\\
	|L_{2}\left(\Psi\right)(r,t)| &\le \int_{0}^{(t-r)_+} (1+\tau)^{-(p-1)} d\tau \nonumber\\
	&\quad\times\int_{0}^{t-r-\tau}
	\frac{\lambda |\Psi(\lambda,\tau)|d\lambda}
	{( \sqrt{2r} \ \mbox{or} \ \sqrt{t-r+\lambda-\tau})
	 \sqrt{t-r-\tau-\lambda}}.\nonumber
\end{align}
\end{lem}

\par
In order to construct our solution in the weighted $L^\infty$ space,
we define the following weighted functions:
\begin{align}
	w_{1}(r,t)&:=\tau_{+}(r,t)^{1/2} \tau_{-}(r,t)^{1/2},\label{weight1}\\
	w_{2}(r,t)&:=\tau_{+}(r,t)^{p_1}  \tau_{-}(r,t)^{p_2}\left(\log 2 \frac{\tau_{+}(r,t)}{\tau_{-}(r,t)}\right)^{-p_3} 
      (\log \tau_{-}(r,t))^{-p_4},   
	\label{weight2}\\
      w_{3}(r,t)&:=\tau_{+}(r,t)^{1/2} \tau_{-}(r,t)^{3/2}\label{weight3},
\end{align}
where we set
\begin{align}
	\tau_+(r,t):=\frac{t+r+2k}{k},
	\quad \tau_-(r,t):=\frac{t-r+2k}{k} \label{dai2.1}
\end{align}
and
\begin{align}
	&p_{1}:=\min\left\{\frac{3p-4}{2}, \frac{1}{2} \right\},
	\quad p_{2}:=\max \left\{ 0, \frac{3p-5}{2} \right\},\nonumber\\
	&p_{3}:=\left\{
	\begin{array}{ll}
		0 & (p\neq5/3),\\
		1 & (p=5/3),
	\end{array}
	\right.
\qquad p_{4}:=\left\{
	\begin{array}{ll}
		0 & (1<p<2),\\
		1 & (p=2).
	\end{array}
	\right.\label{dai2.6}
\end{align}
We remark that $w_2$ can be described as
\begin{align}\label{dai2.8}
	w_{2}(r,t)^{-1}&=
	\left\{
	\begin{array}{lll}
	\tau_{+}(r,t)^{(4-3p)/2}   &
	\d \left(1<p < 5/3 \right),\\
     \d \tau_{+}(r,t)^{-1/2} \log  \left(2 \frac{\tau_{+}(r,t)}{\tau_{-}(r,t)}\right)&
     \d \left(p=5/3 \right),\\   
	\tau_{+}(r,t)^{-1/2}  \tau_{-}(r,t)^{(5-3p)/2} &
	\d \left(5/3 <p < 2 \right),\\
	\tau_{+}(r,t)^{-1/2}  \tau_{-}(r,t)^{-1/2} \log \tau_{-}(r,t)&
	\d \left(p = 2 \right).
\end{array}
	\right.
\end{align}
For these weighted functions, we denote the weighted $L^\infty$ norms of $V$ by 
\[
	\|V\|_{i}:=\sup_{(x,t)\in\R^2\times[0,T)}\{w_{i}(|x|,t)|V(x,t)| \}
\]
where $i=1,2,3$.
\par
Finally, we shall introduce some useful representations for $L$.
It is trivial that $1+\tau \ge (2k+\tau)/(2k)$ is valid for $\tau \ge 0$ and $k\ge1$.
Setting $\tau=(\alpha +\beta)/2 \ge 0$ with $\beta \ge -k$,
we have
\begin{align}
	1+\tau \ge \frac{\alpha+2k}{4k}.\label{dai2.3-1}
\end{align}
Changing the variables by 
\begin{align}
	&\alpha=\tau+\lambda, \quad \beta=\tau-\lambda
	\label{dai2.3}
\end{align}
and extending the domain of $(\alpha,\beta)$-integration,
we obtain from Lemma \ref{lem:int_new_L} and (\ref{dai2.3-1})
\begin{align}
	L_{1}\left(\Psi\right)(r,t) \leq
	\frac{Ck}{\sqrt{r}} \int_{-k}^{t-r} d \beta \int_{t-r}^{t+r}
	\frac{ \{(\alpha+2k)/k\}^{2-p} | \Psi^{*}(\alpha,\beta) |}
	{\sqrt{\alpha-(t-r)}} d \alpha\label{dai2.4}
\end{align}
and
\begin{align}
	L_{1}\left(\Psi\right)(r,t)\leq
	C \sqrt{k} \int_{-k}^{t-r} d \beta \int_{t-r}^{t+r}
	\frac{ \{(\alpha+2k)/k\}^{(3-2p)/2} | \Psi^{*}(\alpha,\beta) |  }
	{\sqrt{\alpha-(t-r)}} d \alpha,\label{dai2.4-1}
\end{align}
where $\Psi^{*}(\alpha,\beta):=\Psi \left((\alpha-\beta)/2,(\alpha+\beta)/2\right)$.\\
Similarly, we get
\begin{align}
	L_{2}\left(\Psi\right)(r,t) \leq 
	\frac{Ck}{\sqrt{r}} \int_{-k}^{t-r} d \beta \int_{-k}^{t-r}
	\frac{ \left\{(\alpha+2k)/k \right\}^{2-p}  | \Psi^{*} (\alpha,\beta) |}
	{ \sqrt{t-r-\alpha}} d\alpha
	\label{dai2.5}
\end{align}
and
\begin{align}
	L_{2}\left(\Psi\right)(r,t) \leq 
	Ck \int_{-k}^{t-r} d \beta \int_{-k}^{t-r}
	\frac{ \left\{(\alpha+2k)/k \right\}^{2-p} | \Psi^{*} (\alpha,\beta) |}
	{ \sqrt{t-r-\alpha} \sqrt{t-r-\beta}} d\alpha.
	\label{dai2.5-1}
\end{align}

\section{Proof of Theorem \ref{T.1.1} and Theorem \ref{T.1.2}}
\par
In this section, we prove Theorem \ref{T.1.1} and Theorem \ref{T.1.2}. 
The proof is based on the classical iteration method in John \cite{John79}.
Lemma \ref{lem3.1} will be used to prove Theorem \ref{T.1.1}, whereas
we prove Theorem \ref{T.1.2} by using Lemma \ref{lm:apriori3}.
First, we prepare the elementary inequalities in Lemma \ref{lem3.4} and Lemma \ref{lem3.5}.
\begin{lem}\label{lem3.4}
Let $a_{1} \in \R$ and $k\ge1$. For $0 \le r \le t+k$, it holds
\begin{align}
	\int_{t-r}^{t+r} \frac{ \left\{ (\alpha+2k)/k \right\}^{a_{1}}}{ \sqrt{\alpha-(t-r)}} d \alpha
	\leq C\sqrt{k} \times \left\{
	\begin{array}{ll}
		\tau_{+}(r,t) ^{a_1+1/2} & (a_1>- 1/2),\\
		\log \left( 2 \d\frac{\tau_{+}(r,t)}{\tau_{-}(r,t)} \right) & (a_1=-1/2),\\
		\tau_{-}(r,t)^{a_1+1/2} & (a_1<-1/2),
	\end{array}
	\right.\nonumber
\end{align}
where $\tau_{+}(r,t)$ and $\tau_{-}(r,t)$ are defined in (\ref{dai2.1}).
\end{lem}
\par\noindent
{\bf Proof.}
For $0 \le r \le t+k$, the integration by parts yields
\begin{align*}
	&\int_{t-r}^{t+r} \frac{ \left\{ (\alpha+2k)/k \right\}^{a_1}}{\sqrt{\alpha-(t-r)}} 
	d\alpha\\
	&\leq 2 \sqrt{2r} \left(\frac{t+r+2k}{k} \right)^{a_1}+\frac{2 | a_1|}{\sqrt{k}}
	\int_{t-r}^{t+r} \left(\frac{\alpha+2k}{k} \right)^{a_1-1/2} d \alpha\\
	&\leq 2\sqrt{2} \sqrt{k} \tau_{+}(r,t)^{a_1+1/2}\\
	& \quad + 2 |a_1| \sqrt{k} \times \left\{
	\begin{array}{ll}
	\d\frac{1 }{a_1+1/2}
	\left(\frac{t+r+2k}{k} \right)^{a_{1}+1/2} & (a_1 > -1/2),\\
	\log\left( \d\frac{t+r+2k}{t-r+2k} \right) & (a_1=-1/2),\\
	\d\frac{1}{|a_1+1/2|}
	\left(\frac{t-r+2k}{k} \right)^{a_{1}+1/2} & (a_1 < -1/2).
	\end{array}
	\right.
\end{align*}
This completes the proof.
\hfill$\Box$
\vskip10pt
\begin{lem}\label{lem3.5}
Let $a_{1} \in \R$ and $k \ge1 $. For $0 \le r \le t+k$, it holds
\begin{align}
	\int_{-k}^{t-r} \frac{\left\{ (\alpha+2k)/k  \right\}^{a_{1}}}{\sqrt{t-r-\alpha}} d \alpha
	\leq C\sqrt{k} \times
	\left\{
	 \begin{array}{ll}
	 \tau_{-}(r,t) ^{ a_{1}+1/2 } & (a_{1}>-1),\\
	 \tau_{-}(r,t)^{-1/2} \log \tau_{-}(r,t) & (a_{1}=-1),\\
	 \tau_{-}(r,t)^{-1/2} & (a_{1}<-1),
	 \end{array}
	\right.\label{dai3.19}
\end{align}
where $\tau_{+}(r,t)$ and $\tau_{-}(r,t)$ are defined in (\ref{dai2.1}).
\end{lem}
\par\noindent
{\bf Proof.}
For $a_{1} \geq 0$, we obtain
\[
\begin{array}{lll}
  \d	\int_{-k}^{t-r} \frac{\left\{(\alpha+2k)/k \right\}^{a_{1}}}{\sqrt{t-r-\alpha}} d \alpha
  	& \d \leq \left(\frac{t-r+2k}{k}\right)^{a_{1}} \int_{-k}^{t-r} \frac{1}{\sqrt{t-r-\alpha}}
  	 d \alpha\nonumber\\
  	 & \d \leq 2 \sqrt{k}  \tau_{-}(r,t) ^{a_{1}+1/2}.
\end{array}
\]
Hence, we obtain (\ref{dai3.19}) for $a_1 \geq 0$.

For $a_1<0$, we show (\ref{dai3.19}). Let $-k \le t-r \le k$, i.e., $k \le t-r+2k \le 3k$.
It follows that
\begin{align}
	\int_{-k}^{t-r} \frac{ \left\{ (\alpha+2k)/k \right\}^{a_1}}{\sqrt{t-r-\alpha}} d \alpha
	& \leq \int_{-k}^{t-r}
	\frac{1}{\sqrt{t-r-\alpha}} d\alpha\nonumber\\
	& \leq   3^{-a_1} \sqrt{k} \tau_{-}(r,t)^{a_{1}+1/2}.\nonumber
\end{align}
We get (\ref{dai3.19}) for $a_1 <0$ and $-k \le t-r \le k$.
\par
Let $t-r \ge k$ which implies $t-r \ge (t-r+2k)/4$. Then, breaking the integral up into 
two pieces, we get
\begin{align}
	&\d\int_{-k}^{t-r} \frac{ \left\{ (\alpha+2k)/k \right\}^{a_1}}
	{\sqrt{t-r-\alpha}} d \alpha\nonumber\\
	&\d = \int_{-k}^{(t-r)/2} \frac{\left\{ (\alpha+2k)/k \right\}^{a_1}}
	{\sqrt{t-r-\alpha}} d \alpha
	+\int_{(t-r)/2}^{t-r} 
	\frac{\left\{ (\alpha+2k)/k \right\}^{a_1}}{\sqrt{t-r-\alpha}} d \alpha\nonumber\\
	&=: J_{1}+J_{2}.\label{dai3.14-1}
\end{align}
It is easy to see that
\begin{align}
	J_{1} &\d \leq \sqrt{2} (t-r)^{-1/2} 
	\int_{-k}^{(t-r)/2} \left(\frac{\alpha+2k}{k}\right)^{a_1} d\alpha\nonumber\\
	&\leq 2\sqrt{2} \sqrt{k} \times
	\left\{
	\begin{array}{ll}
	\d\frac{1}{a_1+1} 
	\tau_{-}(r,t) ^{a_1+1/2} & (a_1> -1),\\
	 \tau_{-}(r,t)^{-1/2} \log \tau_{-}(r,t) & (a_{1}=-1),\\
	\d\frac{1}{|a_{1}+1|} \tau_{-}(r,t)^{-1/2} & (a_{1}<-1).
	\end{array}
	\right.\label{dai3.14-2}
\end{align}
We obtain
\begin{align}
	J_{2} &\d\leq \left\{  \left(\frac{t-r}{2}+2k\right)/k \right\}^{a_1}
	\int_{(t-r)/2}^{t-r} \frac{1}{\sqrt{t-r-\alpha}} d\alpha\nonumber\\
	&\d\leq 2^{-a_1}
	 \left(\frac{t-r+2k}{k}\right)^{a_1} \sqrt{2} (t-r)^{1/2}\nonumber\\
	&\d\leq 2^{-a_1+1/2} \sqrt{k} \tau_{-}(r,t)^{a_1+1/2}.\label{dai3.14-3}
\end{align}
By (\ref{dai3.14-1}), (\ref{dai3.14-2}) and (\ref{dai3.14-3}), we obtain the desired inequality in (\ref{dai3.19}) for $a_{1} <0$ and $t-r \ge k$.
This completes the proof.
\hfill$\Box$
\vskip10pt
The following lemma contains one of the most essential estimates.
\begin{lem}\label{lem3.1}
Let $1< p \le 2$ and $L$ be the linear integral operator defined by {\rm (\ref{L})}.
Assume that $V\in C(\R^2\times[0,T))$ with
{\rm  supp} $V\subset\{(x,t)\in \R^2 \times[0,\infty) : |x|\le t+k\}$.
Then, there exists a positive constant $C_1$ independent of $k$ and $T$ such that 
\begin{equation}
\label{apriori1}
\|L(|V|^p)\|_{1} \le C_1k^2\|V\|_{1}^pD_{1}(T),
\end{equation}
where $D_{1}(T)$ is defined by 
\begin{align}
D_{1}(T):=
\left\{
\begin{array}{lll}
	\d T_{k}^{4-2p} 
	& \mbox{if}\quad 1<p<2,\\
	\d (\log T_{k})^2
	& \mbox{if}\quad p=2
\end{array}
\right.\label{dai3.1-1}
\end{align}
with $T_{k}:=(T+3k)/k.$
\end{lem}
\par\noindent
{\bf Proof.}
In order to show the a-priori estimate (\ref{apriori1}), it is enough to prove
\begin{align}
	L_{j}(w_{1}^{-p}) \leq C k^{2} w_{1}^{-1} D_{1}(T) \quad \mbox{for} \ j=1,2,
	\label{dai3.1}
\end{align}
where $L_{j}$ are defined in Lemma \ref{lem:int_new_L}.
By (\ref{weight1}), (\ref{dai2.1}) and (\ref{dai2.3}), we have
\begin{align}
w_{1}(\lambda,\tau)= \left(\frac{\alpha+2k}{k}\right)^{1/2} 
\left(\frac{\beta+2k}{k}\right)^{1/2}.\label{dai3.2}
\end{align}
We shall prove (\ref{dai3.1}) in the following two cases.
\par
\vskip10pt
\noindent 
Case 1: $4r \geq t+r+2k$.
\par
First, we evaluate $L_{1}$.
We get from (\ref{dai2.4}) and (\ref{dai3.2})
\begin{align}
	L_{1}(w_{1}^{-p}) \leq  \frac{Ck}{\sqrt{r}}
	\int_{t-r}^{t+r} \frac{ \left\{ (\alpha+2k)/k \right\}^{(4-3p)/2}}{\sqrt{\alpha-(t-r)}} 
	d \alpha
	\int_{-k}^{t-r} \left(\frac{\beta+2k}{k}\right)^{-p/2}d \beta.\label{dai3.3}
\end{align}
From Lemma \ref{lem3.4} and (\ref{log-est}), we obtain
\begin{equation}\label{dai3.5}
\begin{array}{llll}
\d	\int_{t-r}^{t+r} \frac{\left\{ (\alpha+2k)/k \right\}^{(4-3p)/2}}
	{ \sqrt{\alpha-(t-r)}} d\alpha\\
	\leq 
C \sqrt{k} \times
	\left\{	
	\begin{array}{ll}
	  	\tau_{+}(r,t)^{(5-3p)/2} & (1<p<5/3),\\
	 	\tau_{+}(r,t)^{1/2} \tau_{-}(r,t)^{-1/2} & (p=5/3),\\
	 	\tau_{-}(r,t)^{(5-3p)/2} & (5/3< p \leq 2).
	 \end{array}
	 \right.
\end{array}
\end{equation}
Here, for the inequality with $p=5/3$, we used the following fact:
\begin{equation}
\label{log-est}
\log s \leq \frac{s^{\delta}}{\delta} \ \mbox{for} \ s \geq 1\ \mbox{and} \ \delta>0.
\end{equation}
The $\beta$-integral is estimated by
\begin{align}
	\int_{-k}^{t-r} \left(\frac{\beta+2k}{k}\right)^{-p/2} d\beta \leq
	Ck \times
	\left\{ 
	\begin{array}{ll}
	\tau_{-}(r,t)^{(2-p)/2} & (1<p<2),\\
	\log \tau_{-}(r,t) & (p=2).	
	\end{array}
	\right.\label{dai3.4}
\end{align}
Therefore, it follows from (\ref{dai3.3}), (\ref{dai3.5}), (\ref{dai3.4}), (\ref{weight1}) and (\ref{dai3.1-1}) that
\[
\begin{array}{lll}
	L_{1}(w_{1}^{-p}) &\leq Ck^2 \tau_{+}(r,t)^{-1/2} \tau_{-}(r,t)^{-1/2}
	\times \left\{
	\begin{array}{ll}
	\tau_{+}(r,t)^{4-2p} & (1<p<2),\\
	\log \tau_{-}(r,t)  & (p=2)
	\end{array}
	\right.\nonumber\\
	&\leq Ck^2 w_{1}(r,t)^{-1} D_{1}(T).
\end{array}
\]
Here, we have used that
\begin{align*}
	\tau_{+}(r,t) \leq \frac{2t+3k}{k} \leq 2T_{k} \quad \mbox{and} \quad T_{k} \geq 3.
\end{align*}
Thus, we have proved (\ref{dai3.1}) with $j=1$ in Case 1.
\par
Next, if $t>r$, we investigate the integral $L_{2}$. From (\ref{dai2.5})
and (\ref{dai3.2}), we get
\begin{align}
         &L_{2}(w_{1}^{-p})(r,t)\nonumber\\
	&\leq \frac{C k}{\sqrt{r}}
	\int_{-k}^{t-r}\frac{\left\{ (\alpha+2k)/k \right\}^{(4-3p)/2}}{\sqrt{t-r-\alpha}} d\alpha
	\int_{-k}^{t-r} \left( \frac{\beta+2k}{k} \right)^{-p/2} d\beta.
	\label{dai3.10}
\end{align}
From Lemma \ref{lem3.5}, we obtain
\begin{align}
	&\int_{-k}^{t-r} \frac{ \left\{(\alpha+2k)/k \right\}^{(4-3p)/2}}
	{\sqrt{t-r-\alpha}} d \alpha\nonumber\\
	&\leq C \sqrt{k} \times
	\left\{
	\begin{array}{ll}
	\tau_{-}(r,t)^{(5-3p)/2} & (1<p<2),\\
	\tau_{-}(r,t)^{-1/2} \log\tau_{-}(r,t) & (p=2).
	\end{array}
	\right.\label{dai3.12}
\end{align}
From (\ref{dai3.10}), (\ref{dai3.12}), (\ref{dai3.4}),
(\ref{weight1}) and (\ref{dai3.1-1}), it follows that
\begin{align}
	L_{2}(w_1^{-p}) &\leq Ck^2 \tau_{+}(r,t)^{-1/2}\times
	\left\{
	\begin{array}{ll}
	\tau_{-}(r,t)^{(7-4p)/2} & (1<p<2),\\
	\tau_{-}(r,t)^{-1/2} \left\{ \log \tau_{-}(r,t) \right\}^2 & (p=2)
	\end{array}
	\right.\nonumber\\
	&\leq  Ck^2 w_{1}(r,t)^{-1}D_{1}(T).\nonumber
\end{align}
Hence, we obtain (\ref{dai3.1}) with $j=2$ in Case 1. 
\par
\vskip10pt
\noindent
Case 2: $4r \le t+r+2k$, i.e., $t+r+2k \le 2(t-r+2k)$.
\par
First, we estimate $L_1$. We have from (\ref{dai2.4-1}) and (\ref{dai3.2})
\begin{align}
	&L_{1}(w_{1}^{-p})\nonumber\\
	&\leq C \sqrt{k}\int_{t-r}^{t+r} \frac{ \left\{ (\alpha+2k)/k \right\}^{3(1-p)/2}}{\sqrt{\alpha-(t-r)}} d \alpha
	\int_{-k}^{t-r} \left( \frac{\beta+2k}{k} \right)^{-p/2} d \beta.
	\label{dai3.8}
\end{align}
We obtain
\begin{align}
	&\int_{t-r}^{t+r} \frac{ \left\{ (\alpha+2k)/k \right\}^{3(1-p)/2}}{\sqrt{\alpha-(t-r)}} d\alpha\nonumber\\
	&\leq \left(\frac{t-r+2k}{k}\right)^{3(1-p)/2} \int_{t-r}^{t+r} 
	\frac{1}{\sqrt{\alpha-(t-r)}} d \alpha\nonumber\\
	&\leq 2 \sqrt{2r} \tau_{-}(r,t)^{3(1-p)/2}\nonumber\\
	&\leq C \sqrt{k} \tau_{+}(r,t)^{(4-3p)/2}.
	\label{dai3.9}
\end{align}
Therefore, it follows from (\ref{dai3.8}), (\ref{dai3.9}), (\ref{dai3.4}),
(\ref{weight1}) and (\ref{dai3.1-1})
\begin{align}
	L_{1}(w_{1}^{-p}) & \leq Ck^2 \tau_{+}(r,t)^{-1/2} \tau_{-}(r,t)^{-1/2}
	\times
	\left\{
	\begin{array}{ll}
	\tau_{+}(r,t)^{4-2p} & (1<p<2),\\
	\log \tau_{-}(r,t) & (p=2)
	\end{array}
	\right.\nonumber\\
	& \leq Ck^2 w_{1}(r,t)^{-1} D_{1}(T).\nonumber
\end{align}
Thus, we obtain (\ref{dai3.1}) with $j=1$ in Case 2.
\par
Next, if $t>r$, we evaluate $L_2$. From (\ref{dai2.5-1}) and (\ref{dai3.2}), we obtain
\begin{align}
	&L_{2}(w_{1}^{-p})(r,t)\nonumber\\
	&\leq C k
	\int_{-k}^{t-r}\frac{\left\{ (\alpha+2k)/k \right\}^{(4-3p)/2}}{\sqrt{t-r-\alpha}} d\alpha
	\int_{-k}^{t-r} \frac{ \left\{(\beta+2k)/k\right\}^{-p/2}}{\sqrt{t-r-\beta}} d\beta.
	\label{dai3.10-1}
\end{align}
From Lemma \ref{lem3.5}, we have
\begin{align}
	&\d\int_{-k}^{t-r} \frac{ \left\{ (\beta+2k)/k \right\}^{-p/2}}{\sqrt{t-r-\beta}} 
	d \beta\nonumber\\
	&\d\leq C \sqrt{k} \times
	\left\{
	\begin{array}{ll}
	\tau_{-}(r,t)^{-(p-1)/2} & (1<p<2),\\
	\tau_{-}(r,t)^{-1/2} \log\tau_{-}(r,t) & (p=2).
	\end{array}
	\right.\label{dai3.12-1}
\end{align}
From (\ref{dai3.10-1}), (\ref{dai3.12}), (\ref{dai3.12-1}), (\ref{weight1}) and (\ref{dai3.1-1}), we have
(\ref{dai3.1}) with $j=2$ in Case 2.
Therefore, the proof of Lemma \ref{lem3.1} is completed.
\hfill$\Box$
\vskip10pt
\par
Finally, we state an a-priori estimate of mixed type.
\begin{lem}
\label{lm:apriori3}
Let $1<p \le 2$ and $L$ be the linear integral operator defined by {\rm (\ref{L})}.
Assume that $V,V_0 \in C(\R^2 \times[0,T))$
with {\rm  supp} $(V,V_0)\subset\{(x,t)\in \R^2 \times[0,\infty) : |x|\le t+k\}$. 
Then, there exists a positive constant $C_2$ independent of $k$ and $T$ such that 
\begin{equation}\label{apriori3}
	\|L(|V_0|^{p-\nu}|V|^{\nu})\|_{2} \le C_2 k^2 \|V_0\|_3^{p-\nu}\| V \|_{2}^{\nu} D_{2,\nu}(T),
\end{equation}
where $\nu=0,p-1,1, p$ and 
\begin{align}
D_{2,\nu}(T)
:=\left\{
\begin{array}{ll}
T_k^{\nu(5-3p)/2+\delta \nu p_3} & (\nu=0, p-1\ \mbox{and}\ 1<p \le 5/3),\\
T_k^{5-3p+\delta \nu p_3} & (\nu=1 \ \mbox{and}\ 1<p \le 5/3),\\
1 & (\nu \le 1 \ \mbox{and} \ 5/3< p\le2),\\
T_k^{\gamma(p,4)/2}& (\nu=p \ \mbox{and} \ 1<p<2),\\
(\log T_k)^3& (\nu=p \ \mbox{and} \ p=2),\\
\end{array}
\right.\label{dai3.21}
\end{align}
where $\delta$ stands for any positive constant and $p_3$ is defined in (\ref{dai2.6}).
\end{lem}

\noindent
{\bf Proof.}
In order to show the a-priori estimate (\ref{apriori3}), it is enough to prove
\begin{equation}
\label{Basic-est2}
L_{j}(w_{3}^{-(p-\nu)} w_{2}^{-\nu}) \leq C k^{2} w_{2}^{-1} D_{2,\nu}(T) \quad \mbox{for} \ j=1,2,
\end{equation}
where $L_j$ are defined in Lemma \ref{lem:int_new_L}. 
For $\delta >0$, from (\ref{weight2}), (\ref{weight3}), (\ref{dai2.1}), (\ref{dai2.3})
and (\ref{log-est}), we have
\begin{align}
	&w_{3}(\lambda,\tau)^{-(p-\nu)} w_{2}(\lambda,\tau)^{-\nu}\nonumber\\
	&\le C  
	\left(\frac{\alpha+2k}{k}\right)^{-(p-\nu)/2-\nu p_1+\delta \nu p_3} 
	\left(\frac{\beta+2k}{k}\right)^{-3(p-\nu)/2-\nu p_2-\delta \nu p_3}\nonumber\\
	&\qquad \times
	\left(\log\frac{\beta+2k}{k}\right)^{\nu p_4}.
	\label{dai3.31}
\end{align}
We shall prove (\ref{Basic-est2})  in the following two cases.
\vskip10pt
\par
\noindent 
Case 1: $4r \geq t+r+2k$.
\par
First, we evaluate $L_1$. From (\ref{dai2.4}) and (\ref{dai3.31}), we get
\begin{align}
	&L_{1}(w_{3}^{-(p-\nu)} w_2^{-\nu})\nonumber\\
	&\leq \frac{Ck }{\sqrt{r}}
	\int_{t-r}^{t+r} \frac{ \{(\alpha+2k)/k\}^{p_5}}
	{\sqrt{\alpha-(t-r)}} d \alpha
	\int_{-k}^{t-r}\left(\frac{\beta+2k}{k}\right)^{p_6}  \left( \log \frac{\beta+2k}{k} \right)^{\nu p_4} d \beta,
	\label{dai3.22}
\end{align}
where
\begin{align}
	p_5&:=\frac{4-3p}{2}+\nu\left(\frac{1}{2}-p_1\right)+\delta \nu p_3,
	\label{dai3.34}\\
	p_6&:=-\frac{3(p-\nu)}{2}-\nu p_2 -\delta \nu p_3.\label{dai3.35}
\end{align}
We have from (\ref{dai3.34}) and (\ref{dai2.6})
\begin{align}
	p_5= \left\{
	\renewcommand{\arraystretch}{1.5}
	\begin{array}{ll}
		\d\frac{4-3p}{2}+\nu \left( \frac{5-3p}{2} \right) & (1< p < 5/3),\\
		\d-\frac{1}{2}+\delta \nu & (p=5/3),\\
		\d\frac{4-3p}{2} & (5/3 < p \le 2).
	\end{array}
	\renewcommand{\arraystretch}{1}
	\right.\label{dai3.23}
\end{align}
From Lemma \ref{lem3.4} and (\ref{dai3.23}), we obtain
\begin{align}
	&\int_{t-r}^{t+r} \frac{\{(\alpha+2k)/k\}^{p_5}}
	{ \sqrt{\alpha-(t-r)}} d\alpha\nonumber\\
	&\leq C \sqrt{k} \times
	\left\{ 
	\begin{array}{ll}
	  	\tau_{+}(r,t)^{(1+\nu)(5-3p)/2} & (1<p<5/3),\\
	 	\log \left(\d 2 \frac{\tau_{+}(r,t)}{\tau_{-}(r,t)}\right) & (\nu=0 \ \mbox{and} \ p=5/3),\\
	 	\tau_{+}(r,t)^{\delta \nu} & (\nu>0 \ \mbox{and} \ p=5/3),\\
	 	\tau_{-}(r,t)^{(5-3p)/2}& (5/3<p \le 2).
	\end{array}
	\right.\label{dai3.24}
\end{align}
We get from (\ref{dai2.6}) and (\ref{dai3.35})
\begin{align}
	p_6=\left\{
	\begin{array}{ll}
		-3(p-\nu)/2-\delta \nu p_3 & (1<p \le 5/3),\\
		\!\!\begin{array}{l}
		-3(p-\nu)/2
		-\nu(3p-5)/2
		\end{array} & (5/3 < p \le 2).
	\end{array}
	\right.\label{dai3.25}
\end{align}
From (\ref{dai3.25}), (\ref{dai2.6}) and (\ref{gamma}), the $\beta$-integral is estimated by
\begin{align}
& \d \int_{-k}^{t-r}\left(\frac{\beta+2k}{k}\right)^{p_6}
\left( \log \frac{\beta+2k}{k} \right)^{\nu p_4} d\beta \nonumber\\
& \d \leq Ck \times
	\left\{
	\begin{array}{ll}
	 	1 & (\nu=0 \ \mbox{or} \ \nu=p-1),\\
		\tau_{-}(r,t)^{(5-3p)/2+p_2} & 
		(\nu=1),\\
		\tau_{-}(r,t)^{1-\delta \nu p_3} & 
		(\nu=p \ \mbox{and} \ 1< p \le 5/3),\\
		\tau_{-}(r,t)^{\gamma(p,4)/2} & 
		(\nu=p \ \mbox{and} \ 5/3<p<2),\\
		\left(\log \tau_{-}(r,t) \right)^{3} & 
		(\nu=p \ \mbox{and} \ p=2).\\
		\end{array}
	\right.\label{dai3.26}
\end{align}
Here, for the inequality with $\nu=p$ and $p=5/3$, we took $0 < \delta \nu <1$.
It follows from (\ref{dai3.22}), (\ref{dai3.24}), (\ref{dai3.26}), (\ref{dai2.6}), (\ref{dai2.8}) and
(\ref{dai3.21}) that
\begin{align}
	L_1(w_3^{-(p-\nu)} w_2^{-\nu})
	\leq & \ C k^2  w_2(r,t)^{-1}\nonumber\\
	&\times
	\left\{
	\begin{array}{ll}
		T_{k}^{\nu(5-3p)/2 + \delta \nu p_3} & (\nu=0, p-1 \ \mbox{and} \ 1 < p \le 5/3),\\
	T_{k}^{5-3p+ \delta \nu p_3} & (\nu=1 \ \mbox{and} \ 1 < p \le 5/3),\\
	1 & (\nu \le 1 \ \mbox{and} \ 5/3 < p \le 2),\\
	T_{k}^{\gamma(p,4)/2} & (\nu=p \ \mbox{and} \ 1 < p < 2),\\
	\left( \log T_{k} \right)^2 & (\nu=p \ \mbox{and} \  p =2)\\
	\end{array}
	\right.\nonumber\\
	\leq & \ Ck^2 w_{2}(r,t)^{-1} D_{2,\nu}(T).\nonumber
\end{align}
We obtain (\ref{Basic-est2}) with $j=1$ in Case 1.
\par
Next, if $t>r$, we investigate $L_{2}$.
From (\ref{dai2.5}) and (\ref{dai3.31}), we get
\begin{align}
	&L_{2}(w_{3}^{-(p-\nu)}w_{2}^{-\nu})(r,t)\nonumber\\
	& \leq   \frac{Ck}{\sqrt{r}}
	\int_{-k}^{t-r}\frac{\{(\alpha+2k)/k\}^{p_5}}
	{\sqrt{t-r-\alpha}} d\alpha
	\int_{-k}^{t-r}\left(\frac{\beta+2k}{k}\right)^{p_6}
	\left( \log \frac{\beta+2k}{k} \right)^{\nu p_4}d\beta.
	\label{L_2-alpha2}
\end{align}
From Lemma \ref{lem3.5} and (\ref{dai3.23}), we have
\begin{align}
	&\int_{-k}^{t-r}\frac{\{(\alpha+2k)/k\}^{p_5} d\alpha}
	{\sqrt{t-r-\alpha}}\nonumber\\
	&\le C\sqrt{k} \times
	\left\{
	\begin{array}{ll}
	 	\tau_{-}(r,t)^{(\nu+1)(5-3p)/2+\delta \nu p_3} & (1<p\le 5/3),\\
		\tau_{-}(r,t)^{(5-3p)/2} & (5/3<p < 2),\\
		\tau_{-}(r,t)^{-1/2} \log \tau_{-}(r,t) & (p=2).
	\end{array}
	\right.\label{dai3.27}
\end{align}
Making use of (\ref{L_2-alpha2}), (\ref{dai3.27}), (\ref{dai3.26}), (\ref{gamma}), (\ref{dai2.6}), (\ref{dai2.8}) and (\ref{dai3.21}), we get
\begin{align}
	&L_2(w_{3}^{-(p-\nu)}w_{2}^{-\nu})(r,t)\nonumber\\
	&\leq Ck^2 \tau_{+}(r,t)^{-1/2}\nonumber\\
	&\qquad\times \left\{
	\begin{array}{ll}
	\tau_{-}(r,t)^{(\nu+1)(5-3p)/2+\delta \nu p_3} & (\nu=0,p-1 \ \mbox{and} 
	\ 1<p \le 5/3),\\
	\tau_{-}(r,t)^{3(5-3p)/2+\delta \nu p_3} & (\nu=1 \ \mbox{and} \ 1<p \le 5/3),\\
	\tau_{-}(r,t)^{(5-3p)/2} & (\nu \le 1 \ \mbox{and} \ 5/3 < p <2),\\
	\tau_{-}(r,t)^{-1/2} \log \tau_{-}(r,t) & (\nu \le 1 \ \mbox{and} \ p=2),\\
	\tau_{-}(r,t)^{(5-3p)/2+\gamma(p,4)/2} & (\nu = p \ \mbox{and} \ 1 < p <2),\\
	\tau_{-}(r,t)^{-1/2} 
	\left( \log \tau_{-}(r,t) \right)^4 & (\nu =p \ \mbox{and} \ p=2)\\
	\end{array}
	\right.\nonumber\\
	&\leq Ck^2 w_{2}(r,t)^{-1}D_{2, \nu}(T).\nonumber
\end{align}
Hence, we obtain (\ref{Basic-est2}) with $j=2$ in Case 1.
\vskip10pt
\par
\noindent
Case 2: $4r \leq t+r+2k$, i.e., $t+r+2k<2(t-r+2k)$.
\par
First, we evaluate $L_1$.
From (\ref{dai2.4-1}) and  (\ref{dai3.31}), we have
\begin{align}
	&L_{1}(w_{3}^{-(p-\nu)}w_{2}^{-\nu})\nonumber\\
	& \leq C\sqrt{k}
	\int_{t-r}^{t+r} \frac{\{(\alpha+2k)/k\}^{3(1-p)/2+\nu(1/2-p_1)+\delta \nu p_3}}{\sqrt{\alpha-(t-r)}} 
	d \alpha\nonumber\\
	&\qquad\times \int_{-k}^{t-r} \left(\frac{\beta+2k}{k}\right)^{p_6}
	\left( \log \frac{\beta+2k}{k} \right)^{\nu p_4} d \beta
	\label{dai3.28}.
\end{align}
From (\ref{dai3.34}), we obtain
\begin{align}
	&\int_{t-r}^{t+r} \frac{\{(\alpha+2k)/k\}^{3(1-p)/2+\nu(1/2-p_1)+\delta \nu p_3}}
	{\sqrt{\alpha-(t-r)}} d\alpha\nonumber\\
	&\leq \tau_{-}(r,t)^{3(1-p)/2}  \tau_{+}(r,t)^{\nu(1/2-p_1)+\delta \nu p_3}\int_{t-r}^{t+r} 
	\frac{1}{\sqrt{\alpha-(t-r)}} d \alpha\nonumber\\
	&\leq C \sqrt{k} \tau_{+}(r,t)^{p_5}\label{dai3.29}.
\end{align}
Therefore, it follows from (\ref{dai3.28}), (\ref{dai3.29}), (\ref{dai3.23}), (\ref{dai3.26}),
(\ref{dai2.8}) and (\ref{dai3.21}) that
\begin{align}
	L_{1}(w_{3}^{-(p-\nu)}w_{2}^{-\nu})  \leq
	& \ Ck^2 \tau_{+}(r,t)^{(4-3p)/2}\nonumber\\
	&\times
	\left\{
	\begin{array}{ll}
		\tau_{+}(r,t)^{\nu(5-3p)/2 + \delta \nu p_3} & (\nu=0, p-1 \ \mbox{and} \ 1<p \le 5/3),\\
		\tau_{+}(r,t)^{5-3p+ \delta \nu p_3} & (\nu=1 \ \mbox{and} \ 1<p\le5/3),\\
		1 & (\nu \le 1 \ \mbox{and} \ 5/3 < p \le 2),\\
		\tau_{-}(r,t)^{\gamma(p,4)/2} & (\nu =p \ \mbox{and} \ 1<p<2),\\
	 	\left( \log \tau_{-}(r,t) \right)^3& (\nu =p \ \mbox{and} \ p=2)\\
	\end{array}
	\right.\nonumber\\
	\leq & \  Ck^2 w_{2}(r,t)^{-1}D_{2,\nu}(T).\nonumber
\end{align}
Thus, the proof of (\ref{Basic-est2}) with $j=1$ in Case 2 is finished.
\par
Next, if $t>r$, we investigate $L_2$. From (\ref{dai2.5-1}), (\ref{dai3.31}), (\ref{dai3.34})
and (\ref{dai3.35}), we get
\begin{align}
	L_{2}(w_{3}^{-(p-\nu)}w_{2}^{-\nu})(r,t) \leq Ck  
	&\int_{-k}^{t-r}\frac{\{(\alpha+2k)/k\}^{p_5}}
	{\sqrt{t-r-\alpha}} d\alpha\nonumber\\
	& \times \int_{-k}^{t-r} 
	\frac{\{(\beta+2k)/k\}^{p_6}}{\sqrt{t-r-\beta}} 
	\left( \log \frac{\beta+2k}{k} \right)^{\nu p_4}
	d\beta
	\label{L_2-beta2}.
\end{align}
From Lemma \ref{lem3.5} and (\ref{dai3.25}), we have
\begin{align}
	 &\int_{-k}^{t-r} \frac{\{(\beta+2k)/k\}^{p_6}}
	{\sqrt{t-r-\beta}} 
	\left( \log \frac{\beta+2k}{k} \right)^{\nu p_4}
	d\beta\nonumber\\
	&\le C\sqrt{k} \times
	\left\{
	\begin{array}{ll}
		\tau_{-}(r,t)^{-1/2} & 
		(\nu=0 \ \mbox{or} \ \nu=p-1),\\
		\tau_{-}(r,t)^{-p_1} & (\nu=1),\\
		\tau_{-}(r,t)^{1/2} & 
		(\nu=p \ \mbox{and} \ 1<p<5/3),\\
		\tau_{-}(r,t)^{-\delta \nu+1/2} & 
		(\nu=p \ \mbox{and} \ p=5/3),\\
		\tau_{-}(r,t)^{-1/2+\gamma(p,4)} & 
		(\nu=p \ \mbox{and} \ 5/3<p<2),\\
		\tau_{-}(r,t)^{-1/2} \log \tau_{-}(r,t) & 
		(\nu=p \ \mbox{and} \ p=2).
	\end{array}
	\right.\label{dai3.33}
\end{align}
Here, for the inequality with $\nu=p$ and $p=5/3$, we took $0< \delta \nu <1$.
Thus, we obtain (\ref{Basic-est2}) with $j=2$ in Case 2 by (\ref{L_2-beta2}), (\ref{dai3.27}), (\ref{dai3.33}), (\ref{dai2.8}) and (\ref{dai3.21}).
Therefore, the proof of Lemma \ref{lm:apriori3} is completed.
\hfill$\Box$
\vskip10pt

In the following, we prove Theorem \ref{T.1.1} and Theorem \ref{T.1.2}.
We remark that it is possible to construct a classical solution if $p=p_S(4)=2$. However, its construction is almost the same as for $C^1$ solution. Therefore, we shall omit the proof. 

\par\noindent
{\bf Proof of Theorem \ref{T.1.1}.}
We define
\begin{align*}
	X: = \left\{
	u(x, t)  \ \Biggl| \
\begin{array}{ll}
	D_x^{\alpha} u(x, t) \in C(\R^2 \times [0, T)),\\
	\| D_x^{\alpha} u\|_1 < \infty 	\quad (|\alpha| \leq 1), \\
	u(x, t) = 0 \quad (|x| \geq t+k) \\
\end{array}
\right\},
\end{align*} 
where $D_x^{\alpha} = D_1^{\alpha_1} D_2^{\alpha_2}$ $(\alpha = (\alpha_1, \alpha_2))$ and  $D_k = \p/\p x_k$ $(k = 1, 2)$.
We can verify easily that $X$ is complete with respect to the norm
\begin{align*}
	\| u \|_{X} =\sum_{| \alpha | \leq 1} \| D_{x}^{\alpha}  u\|_1.
\end{align*}
Using the iteration method, we shall construct a solution of (\ref{IVP2}). We define the sequence of functions $\{u_j\}$ by
\begin{align*}
	u_0 = u^0, \quad u_{j+1} = u^0 + L(|u_{j}|^p) \quad \mbox{for} \quad j \geq 0,
\end{align*}
where $u^{0}$ is defined in (\ref{dai2.7.1}).
It follows from Lemma 1 in \cite{G81}, p.236 that $u^0$ satisfies
\begin{align*}
	|D_{x}^{\alpha} u^0(x,t)| \leq C(f,g) \e (t+r+2k)^{-1/2}
	(t-r+2k)^{-1/2}
\end{align*}
for $|\alpha| \leq 1$, where the positive constant $C(f,g)$ depends on $D^{\alpha}_{x} g$
and $D^{\beta}_{x} f$ \ $(|\beta| \leq 2)$. Hence, we find
\begin{align}
	\| D^{\alpha}_{x} u^0 \|_1 \leq C(f,g) k^{-1} \e.\label{dait10}
\end{align}
As in \cite{John79}, p.258, we see from Lemma \ref{lem3.1} and (\ref{dait10}) that if $\e$ satisfies
\begin{align}
	C_1k^2 D_{1}(T) \e ^{p-1} \| u_L \|^{p-1}_{1} \leq \frac{1}{p2^p},\label{dait8}
\end{align}
then $\{ u_{j} \}$ is a Cauchy sequence in
$X$. Since $X$ is complete, there exists a function $u \in X$ such that $\{ D_{x}^{\alpha} u_j \}$
converges uniformly to $D_{x}^{\alpha} u$ as $j \to \infty$. Clearly $u$ satisfies
 (\ref{dai2.7}) with $F(x,t)=|u(x,t)|^p$. 
In view of (\ref{dai2.7}) and (\ref{L}), we note that
$\partial u/ \partial t$ can be expressed in terms of $D^{\alpha}_x u$ ($|\alpha| \le 1$). 
From the continuity of  $D^{\alpha}_x u$, the continuity of $\partial u/ \partial t$
is also valid. 
Thus, from (\ref{dait10}) and (\ref{dait8}), Theorem \ref{T.1.1} is proved by taking $\e$ is small.
\hfill$\Box$\\

\par\noindent
{\bf Proof of Theorem \ref{T.1.2}.}
We consider the following integral equation:
\begin{align}
	U=L(|u^0+U|^p), \label{dait11}
\end{align}
where $u^{0}$ is defined in (\ref{dai2.7.1}).
Suppose that we obtain a solution $U=U(x,t)$ of (\ref{dait11}). Then, putting $u=U+u^0$, we get the solution
of (\ref{dai2.7}) with $F(x,t)=|u(x,t)|^p$, and its maximal existence time 
is the same as that of $U$. Thus, we have reduced the problem
to the analysis of (\ref{dait11}).
Let $Y$ be the norm space defined by
\begin{align*}
	Y := \left\{
	U(x, t) \ \Biggl| \
\begin{array}{ll}
	D_x^{\alpha} U(x, t) \in C(\R^2 \times [0, T)), \\
	\| D_x^{\alpha} U\|_2 < \infty 	\ (|\alpha| \leq 1), \!\! \\
	U(x, t) = 0 \quad (|x| \geq t+k) \\
\end{array}
\right\},
\end{align*}
which is equipped with the norm
\begin{align*}
	\| U \|_{Y}= \sum_{ |\alpha| \leq 1 } \| D_{x}^{\alpha} U \|_2. 
\end{align*}
We shall construct a solution of the integral equation (\ref{dait11}) in $Y$. We define the sequence of functions $\{ U_{j} \}$ by
\begin{equation}
\label{seq-U_j}
	U_{0}=0,\quad U_{j+1}=L(|u^0+U_{j}|^p) \quad (j=0, 1, 2, \cdots).
\end{equation}
From Lemma \ref{lem:decay_est_v} and the condition $\int_{\R^2}(f+g)(x) dx=0$, we see that there exists a positive constant $C_0$ such that
\begin{equation}
\label{u^0_3}
	\| D_{x}^{\alpha} u^0 \|_{3} \leq C_{0} \e \quad (|\alpha| \leq 1).
\end{equation}
We put
\begin{equation}
\label{C_3}
	C_{3}:=(2^{2p+2}p)^{p/(p-1)}
	\max \left\{ C_{2} k^2 M_{0}^{p-1}, (C_{2} k^2 C_{0}^{p-1})^p,
	(C_{2} k^2 M_{0}^{p-2} C_{0})^{ p/(p-1)} \right\}
\end{equation}
and
\begin{equation}
\label{M_0-def}
	M_{0}:=2^p p k^2 C_{0}^p C_{2},
\end{equation}
where $C_2$ is the constant given in Lemma 
\ref{lm:apriori3}.
We take $\e$ and $T$ such that 
\begin{align}
	C_{3} \e^{p(p-1)} D_{2,p}(T) \leq 1\label{dait9}.
\end{align}
\par
Similarly to the proof of Theorem 1 in \cite{IKTW}, we shall obtain
\begin{equation}
\label{ind-1}
\|U_{j}\|_2\le2M_0\e^p
\end{equation}
by induction. For $j=0$, (\ref{ind-1}) holds. 
Assume that $\| U_{j}\|_2 \leq 2M_{0} \e^{p}$ for some $j$.
From (\ref{seq-U_j}), (\ref{apriori3}) with $\nu=0$ and $\nu=p$, (\ref{u^0_3}), (\ref{C_3}), 
(\ref{M_0-def}) and (\ref{ind-1}), we have 
\begin{align}
\|U_{j+1}\|_2&\le2^{p-1}\{\|L(|u^0|^p)\|_2+\|L(|U_{j}|^p)\|_2\}\nonumber\\
&\le 2^{p-1}C_2k^2\{\|u^0\|_3^pD_{2,0}(T)
+\|U_j\|_2^pD_{2,p}(T)\}\nonumber\\
&\le 2^{p-1}C_2k^2\{(C_0\e)^p+
(2M_0\e^p)^pD_{2,p}(T)\}\nonumber\\
&\le M_0\e^p+M_0C_3\e^{p^2}D_{2,p}(T)\label{ind-1-est}.
\end{align}
Thus, from (\ref{ind-1-est}), we obtain (\ref{ind-1}) 
under the conditions (\ref{dait9}).
\par
Next, we shall estimate the differences of $\{U_j\}$. 
From (\ref{seq-U_j}), we obtain
\begin{align}
	\| U_{j+1} -U_{j} \|_{2} &\leq 2^{p-1} p
	\{2\| L(|u^{0}|^{p-1} | U_{j}-U_{j-1}| )\|_{2}\nonumber\\
	&+\| L((|U_{j}|^{p-1}+|U_{j-1}|^{p-1}) 
	|U_{j}-U_{j-1}|)\|_{2}\}.\label{diff-est-1}
\end{align}
Then, in view of the definitions of $D_{2,1}(T)$ and $D_{2,p}(T)$ in 
(\ref{dai3.21}), the estimate $D_{2,1}(T)\le D_{2,p}(T)^{1/p}$ holds
because of $5-3p+\nu p_3\le \gamma(p,4)/2p$ and $T_k\ge3$. 
Here, when $p=5/3$, we take $\delta>0$ such that $0<\delta<1/p$ 
in Lemma \ref{lm:apriori3}. 
Thus, from (\ref{apriori3}) with $\nu=1$ and (\ref{u^0_3}), we obtain
\begin{align}
\| L(|u^0|^{p-1} |U_{j}-U_{j-1}|)\|_{2}
&\leq C_{2}k^2\| u^{0} \|_{3}^{p-1} \| U_{j} - U_{j-1} \|_{2}D_{2,1}(T)\nonumber\\
&\leq C_{2}k^2(C_0\e)^{p-1} D_{2,p}(T)^{1/p}\| U_{j} - U_{j-1} \|_{2}\label{diff-est-2}.
\end{align}
We also get from (\ref{apriori3}) with $\nu=p$ and (\ref{ind-1}) that
\begin{align}
	&\| L((|U_{j}|^{p-1} +|U_{j-1}|^{p-1})
	|U_{j}-U_{j-1}|)\|_{2}\nonumber\\
	&\leq C_{2}k^2
	(\|  U_{j} \|_{2}^{p-1}+\|  U_{j-1} \|_{2}^{p-1})
	\| U_{j}-U_{j-1} \|_{2}D_{2,p}(T) \nonumber\\
	&\leq 2C_{2}k^2(2M_{0} \e^p)^{p-1}D_{2,p}(T)
	\| U_{j} - U_{j-1} \|_{2}\label{diff-est-3}.
\end{align}
Hence, we obtain from (\ref{diff-est-1}), (\ref{diff-est-2}), 
(\ref{diff-est-3}) and (\ref{C_3}) that 
\begin{align}
\| U_{j+1} - U_{j} \|_{2} 
&\leq \frac{1}{4}\{(C_3\e^{p(p-1)}D_{2,p}(T))^{1/p}
+C_{3}\e^{p(p-1)}D_{2,p}(T)\}\| U_{j}-U_{j-1}\|_{2}
\nonumber\\
&\le \d
\frac{1}{2}\| U_{j}-U_{j-1}\|_{2}\label{diff-est-4}
\end{align}
provided (\ref{dait9}) holds. 

Similarly to the proof of (\ref{ind-1}) and (\ref{diff-est-4}), if we 
assume that (\ref{dait9}) holds, then we obtain the following estimates: 
\begin{equation}
\label{ind-2}
\|D_iU_{j}\|_2\le2M_0\e^p,
\end{equation}
\begin{equation}
\label{diff-est-5}
\| D_{i}(U_{j+1}-U_{j}) \|_{2}
\leq C_4(j+1)2^{-j(p-1)},
\end{equation}
where $C_4$ is a positive constant independent of $j$.
We remark that in order to show (\ref{ind-2}) and (\ref{diff-est-5}), we also use the 
estimates (\ref{apriori3}) with $\nu=p-1$ and $D_{2,p-1}(T) \leq D_{2,p}(T)^{(p-1)/(p+1)}$.
For the actual proof, see the inequalities (4.15) and (4.25) 
in \cite{IKTW} which correspond the estimates (\ref{ind-2}) and (\ref{diff-est-5}) respectively. 
Then, from (\ref{diff-est-4}) and (\ref{diff-est-5}), 
we see that $\{ U_j \}$ is a Cauchy sequence in $Y$ provided that (\ref{dait9}) holds. 
We can verify easily that $Y$ is complete. Hence, there exists a function $U$ such that $\{U_j\}$ converges to $U$ in $Y$. Therefore, $U$ satisfies the integral equation (\ref{dait11}).\\
Let us fix $\e_{0}$ as 
\begin{equation}
\label{dai4-21}
	C_{3} \e_{0}^{p(p-1)}  \leq
	\left\{
	\begin{array}{ll}
		6^{- \gamma(p,4)/2} & (1<p<2),\\
		(\log 6)^{-3} & (p=2).
	\end{array}
	\right.
\end{equation}
For $0<\e\le\e_0$, if we assume that
\begin{align}
	C_{3} \e^{p(p-1)}\leq
	\left\{
	\begin{array}{ll}
	 \left(\d\frac{2T}{k} \right)^{-\gamma(p,4)/2} &(1<p<2),\\ 
	\left( \log\d\frac{2T}{k} \right)^{-3} & (p=2),
	\end{array}
	\right.\label{dai4-20}
\end{align}
then (\ref{dait9}) holds. In fact, when $T\le 3k$,
(\ref{dait9}) follows from (\ref{dai4-21}).
When $T>3k$, (\ref{dait9}) 
follows from (\ref{dai4-20}).
Hence, Theorem \ref{T.1.2} follows immediately 
from (\ref{dai4-20}). This completes the proof.
\hfill$\Box$

\section{Proof of Theorem \ref{T.2.1} and Theorem \ref{T.2.2}}
In this section, we prove Theorem \ref{T.2.1} and Theorem \ref{T.2.2}.
For the sub-critical case, we use an improved version of Kato's lemma on ordinary differential
inequality which was introduced by Takamura \cite{Takamura15}.
For the critical case, we apply the slicing iteration method which was introduced
by Agemi, Kurokawa and Takamura \cite{AKT00}.
From now on, let $u\in C^1(\R^2\times[0, T))$ be the solution of the integral equation associated with (\ref{IVP2}).

\subsection{Proof of Theorem \ref{T.2.1}}
We divide the proof of Theorem \ref{T.2.1} into two cases,
$1<p<2$ and $p=2$. First, we shall handle the sub-critical case.\\
\par\noindent
{\bf Proof of Theorem \ref{T.2.1} with \mbox{\boldmath{$1<p<2$}}.}
We shall follow the arguments in Section 4 of Takamura \cite{Takamura15}.
In order to obtain the estimates in Theorem \ref{T.2.1},
we shall take a look on the ordinary differential inequality for
\begin{align*}
	F(t):=\int_{\R^2}u(x,t)dx.
\end{align*}
(\ref{IVP2}) with $\mu=2$ and (\ref{supp_u}) imply that
\begin{align}\label{dai4.1}
	F''(t)=\frac{1}{(1+t)^{p-1}}\int_{\R^2}|u(x,t)|^pdx
	\quad\mbox{for}\ t\ge0.
\end{align}
Hence, the H\"older's inequality and (\ref{supp_u}) yield that
\begin{align}\label{dai4.2}
	F''(t)\ge{\pi}^{-(p-1)}(t+k)^{-3(p-1)}|F(t)|^p
	\quad\mbox{for}\ t\ge0.
\end{align}
Due to the assumption on the initial data in Theorem \ref{T.2.1},
$f(x) \equiv 0$, $g(x) \geq 0$ $(\not\equiv0)$, we have
\begin{align}
	F(0)=0,\quad F'(0)>0.\label{dai4.3}
\end{align}
It follows from (4.3) in \cite{Takamura15} that
\begin{align}
	u(x,t) \geq  \frac{\| g \|_{L^{1}(\R^2)}}{2\sqrt{2} \pi \sqrt{t+k} \sqrt{t-|x|+k}} 
	 \ \e \quad \mbox{for} \ k \leq |x| \leq t-k.\label{dai4.5}
\end{align}
From (\ref{dai4.1}), it follows that
\begin{align}
	F''(t)  \geq \frac{1}{(1+t)^{p-1}} \int_{k \leq |x| \leq t-k} |u(x,t)|^p dx \quad  \mbox{for} \ t \geq 2k.\nonumber
\end{align}
Plugging (\ref{dai4.5}) into the right-hand side of this inequality, we have that
\begin{align}
	F''(t) & \geq \left( \frac{\| g \|_{L^{1}(\R^2)}}{2 \sqrt{2} \pi (t+k)^{3/2-1/p}}
	\e \right)^{p} \int_{k \leq |x| \leq t-k} \frac{1}{(t-|x|+k)^{p/2}} dx\nonumber\\
	& \quad=\frac{2 \pi \| g \|_{L^{1}(\R^2)}^p}{(2 \sqrt{2} \pi)^p (t+k)^{3p/2-1}}
	\e^p  \int_{k}^{t-k} \frac{r}{(t-r+k)^{p/2}} dr.\label{dai4.7}
\end{align}
We evaluate the integral of the last term in (\ref{dai4.7}). For $t \geq 3k$, we obtain
\begin{align}
	\int_{k}^{t-k} \frac{r}{(t-r+k)^{p/2}} dr &\geq \frac{1}{2t^{p/2}}
	\{(t-k)^2-k^2\}\nonumber\\
	&\geq \frac{1}{6} t^{2-p/2}.\label{dai4.14}
\end{align}
From (\ref{dai4.7}) and (\ref{dai4.14}), we obtain
\begin{align*}
	F''(t) \geq \frac{\| g \|_{L^1(\R^2)}^{p}}{3 \cdot 2^{3p-1} \pi^{p-1}}
	 \, \e^p \, t^{3-2p} \quad \mbox{for} \ t \geq 3k.
\end{align*}
Integrating this inequality in $[3k,t]$, we get from (\ref{dai4.3})
\begin{align*}
	F'(t) > \frac{\| g \|_{L^{1}(\R^2)}^p(1-(3/4)^{4-2p})}
	{3(4-2p)2^{3p-1} \pi^{p-1}} \, \e^p \, t^{4-2p}
	\quad \mbox{for} \ t \geq 4k.
\end{align*}
Hence, we obtain from (\ref{dai4.3})
\begin{align}
	F(t) > D_{1} \e^p t^{5-2p} \quad \mbox{for} \ t \geq 5k,
	\label{dai4.30}
\end{align}
where
\begin{align*}
	D_{1}:=\frac{\| g \|_{L^{1}(\R^2)}^p(1-(3/4)^{4-2p}) (1-(4/5)^{5-2p})}
	{3(4-2p)(5-2p)2^{3p-1} \pi^{p-1}}>0.
\end{align*}
In the sub-critical case, the following basic lemma is useful.
\begin{lem}[Takamura \cite{Takamura15}]\label{lem4.1}
Let $p>1, a>0$ and $q>0$ satisfy
\begin{equation}\label{dai4.8}
	M:=\frac{p-1}{2}a-\frac{q}{2}+1>0.
\end{equation}
Assume that $F\in C^2([0,T))$ satisfies
\begin{align}
	& F(t)\ge  At^a  \quad \mbox{for} \ \ t\ge T_0,\nonumber\\
	& F''(t)\ge  B(t+k)^{-q}|F(t)|^p \quad  \mbox{for} \ \ t\ge0,\nonumber\\
	& F(0) \geq 0,\ F'(0) > 0,\label{dai4.11}
\end{align}
where $A,B,k,T_0$ are positive constants.
Then, there exists a positive constant $D_0=D_0(p,a,q,B)$ such that
\begin{equation}\nonumber
	T<2^{2/M}T_1
\end{equation}
holds provided
\begin{equation}\label{dai4.13}
	T_1:=\max\left\{T_0,\frac{F(0)}{F'(0)},k\right\}\ge D_0 A^{-(p-1)/(2M)}.
\end{equation}
\end{lem}
\noindent
This is exactly Lemma 2.1 in \cite{Takamura15}, so that we shall omit the proof here.
\par
According to (\ref{dai4.2}), (\ref{dai4.3}) and (\ref{dai4.30}), we are in a position to apply our situation to Lemma \ref{lem4.1} with
\begin{align*}
	A=D_{1} \e^p,\quad B=\pi^{1-p},\quad a=5-2p,\quad q=3(p-1)
\end{align*}
which imply that (\ref{dai4.8}) yields
\begin{align*}
	M=\frac{p-1}{2}(5-2p)-\frac{3(p-1)}{2}+1=p(2-p) >0.
\end{align*}
If we set
\begin{align}
	T_{0}:= D_{0} A^{-(p-1)/(2M)}
	= D_{0} D_{1}^{-(p-1)/(p(4-2p))} 
	\e^{-(p-1)/(4-2p)},\label{dai4.15}
\end{align}
we find that there is an $\e_{0}=\e_{0}(g,p,k)$ such that
\begin{align*}
	T_{0} \geq \max\left\{ \frac{F(0)}{F'(0)}, 5k \right\}=5k
	\quad \mbox{for} \ 0<\e \leq \e_{0}.
\end{align*}
This means that $T_1=T_0$ in (\ref{dai4.13}). Therefore, from (\ref{dai4.15}), the conclusion of Lemma \ref{lem4.1} implies 
\begin{align*}
	T < 2^{2/M} T_1=D_{2}\e ^{-(p-1)/(4-2p)},
\end{align*}
where
\begin{align*}
	D_{2}:=2^{2/M} D_{0} D_{1}^{-(p-1)/(p(4-2p))}>0.
\end{align*}
The proof of Theorem  \ref{T.2.1} with $1<p<2$ is now completed.
\hfill$\Box$
\vskip10pt

\par\noindent
{\bf Proof of Theorem \ref{T.2.1} with \mbox{\boldmath{$p=2$}}.}
\par
Let $\bar{v}$ be the spherical mean of $v \in C^{0}(\R^2 \times [0,\infty))$
with radius $r$;
\begin{align}
	\bar{v}(r,t):=\frac{1}{2 \pi} \int_{|\omega|=1} v(r \omega, t) d S_{\omega}.
	\nonumber
\end{align}
We get the following inequality (for the proof, see \cite{AT92}, p.529):
\begin{align}
	\bar{u}(r,t) &\geq \overline{u^0}(r,t)+\frac{2}{\pi} \int_{0}^{t-r} d\tau(1+\tau)^{-1}
	\int_{0}^{t-\tau-r} 
	\lambda |\bar{u}|^2(\lambda,\tau) d\lambda\nonumber\\
	&\quad\times \int_{|\lambda-r|}^{\lambda+r}
	\frac{\rho d\rho}{\sqrt{h(\lambda,\rho;r)((t-\tau)^2-\rho^2)}},
	\label{dai4.61}
\end{align}
where
\begin{align}
	h(\lambda,\rho;r):=(\rho^2-(\lambda-r)^2)((\lambda+r)^2-\rho^2).
	\nonumber
\end{align}
Since $\int_{\R^2} \{ f(x)+g(x) \} dx >0$, from Lemma \ref{lem4.5}, there exist  positive constants $E_0$ and $K$ such that
\begin{align}
	u_{L}(x,t) \geq \frac{E_0}{\sqrt{(t+r)(t-r)}}\nonumber
\end{align}
for $t-r\geq K \geq 1$. Making use of the positivity of the second term of right-hand side in (\ref{dai4.61}),
we get
\begin{align}
	\bar{u}(r,t) &\geq \overline{u^0}(r,t)
	\geq \frac{E_{0} \e }{\sqrt{(t+r)(t-r)}}
	\quad \mbox{for} \ t-r \geq K.\label{dai4.64}
\end{align}
We define
\begin{align}
	\Sigma_{j}:=\left\{ (r,t) \ | \ t-r \geq K l_j\right\},\quad
	\Sigma_{\infty}:=\left\{ (r,t) \ | \ t-r \geq 2K \right\},\nonumber
\end{align}
where
\begin{align*}
	l_j:=\sum_{k=0}^{j}2^{-k}\quad (j=0,1,2,\cdots).
\end{align*}
For $(r,t) \in \Sigma_0$, it follows from (\ref{dai4.61}) and (\ref{dai4.64}) that
\begin{align}
	\bar{u}(r,t) &\geq \frac{2}{\pi} \int_{0}^{t-r} d \tau (1+\tau)^{-1}
	\int_{0}^{t-\tau-r} d\lambda \ \lambda |\bar{u}|^2(\lambda,\tau) 
	\nonumber\\
	&\quad\times \frac{1}{\sqrt{(t-r-(\tau+\lambda))(t+r-(\tau-\lambda))}}\nonumber\\
	&\quad\times \int_{|\lambda-r|}^{\lambda+r} \frac{\rho}{\sqrt{h(\lambda,\rho;r)}} 
	 d\rho.\label{dai4.65}
\end{align}
For $(r,t) \in \Sigma_0$, we introduce
\begin{align}
	Q_{j}(r,t):=\left\{ 
	(\lambda,\tau) \in [0,\infty)^2 \ | \ Kl_{j} \leq \tau-\lambda,
	\ 0 \leq \lambda \leq t-r-\tau
	\right\}.\nonumber
\end{align}
Since
\begin{align}
	\int_{|\lambda-r|}^{\lambda+r} \frac{\rho}{\sqrt{h(\lambda,\rho;r)}}d\rho
	=B\left(\frac{1}{2}, \frac{1}{2}\right)=\frac{\pi}{2},\nonumber
\end{align}
we get from (\ref{dai4.65})
\begin{align}
	\bar{u}(r,t) 
	\geq \frac{1}{\sqrt{(t+r)(t-r)}}
	\iint_{Q_{0}(r,t)} (1+\tau)^{-1} \lambda |\bar{u}|^2 d\lambda d\tau
	\quad \mbox{in} \ \Sigma_0.
	\label{dai4.67}
\end{align}
For $(r,t) \in \Sigma_j$, we have
\begin{align}
	Q_{j}(r,t) \subset Q_{0}(r,t) \quad \mbox{and} \quad Q_{j}(r,t) \subset \Sigma_j.
	\label{dai4.67-1}
\end{align}
Since $\Sigma_j \subset \Sigma_0$, it follows from (\ref{dai4.67}) and
 (\ref{dai4.67-1})
that
\begin{align}
	\bar{u}(r,t) \geq \frac{1}{\sqrt{(t+r)(t-r)}}
	\iint_{Q_{j}(r,t)} (1+\tau)^{-1} \lambda |\bar{u}|^2(\lambda, \tau)
	d\lambda d\tau \quad \mbox{in} \ \Sigma_j.\label{dai4.67-2}
\end{align}
\par
By using the induction argument, we will show
\begin{align}
	\bar{u}(r,t) \geq \frac{d_{j}}{\sqrt{(t+r)(t-r)}} \log^{a_j}
	\left( \frac{t-r}{Kl_j} \right) \quad \mbox{in} \ \Sigma_j,\label{dai4.68}
\end{align}
where
\begin{align}
	&a_{0}=0,\quad a_{j+1}=2a_{j}+2,\label{dai4.69}\\
	&d_0=E_0 \e,\quad d_{j+1}=\frac{d_{j}^2}{2^{3j+9}}.
	\label{dai4.70}
\end{align}
From (\ref{dai4.64}), it holds (\ref{dai4.68}) with $j=0$. We assume 
that (\ref{dai4.68}) holds for one natural number $j$ and $(r,t) \in \Sigma_{j+1}$.
Substituting (\ref{dai4.68})
 into (\ref{dai4.67-2}) and changing the variables by (\ref{dai2.3}), we get
\begin{align}
	&\sqrt{(t+r)(t-r)} \bar{u}(r,t)\nonumber\\
	&\geq \frac{d_j^2}{2}\int_{K l_j}^{t-r} d\alpha
	\int_{K l_j}^{\alpha} \left(1+\frac{\alpha+\beta}{2} \right)^{-1}
	\left(\frac{\alpha-\beta}{2}\right) \alpha^{-1} \beta^{-1}
	\log^{2a_j} \left( \frac{\beta}{K l_j} \right)
	d\beta\nonumber \\
	&\geq \frac{d_j^2}{8(2a_j+1)} \int_{K l_j}^{t-r} d\alpha \alpha^{-2}
	\int_{K l_j}^{\alpha} (\alpha-\beta) \frac{d}{d\beta}
	\left\{  \log^{2a_{j}+1} \left(  \frac{\beta}{K l_j} \right)\right\} d\beta
	\nonumber\\
	&=\frac{d_j^2}{8(2a_j+1)} \int_{K l_j}^{t-r} d\alpha \alpha^{-2}
	\log^{2a_{j}+1} \left(  \frac{\beta}{K l_j} \right)  d\beta
	\nonumber
\end{align}
and, then,
\begin{align}
	&\sqrt{(t+r)(t-r)} \bar{u}(r,t)\nonumber\\
	&\geq\frac{d_j^2}{8(2a_j+1)} \int_{K l_{j+1}}^{t-r} d\alpha \alpha^{-2}
	\int_{\alpha l_{j}/l_{j+1}}^{\alpha} 
	\log^{2a_{j}+1} \left(  \frac{\beta}{K l_j} \right) d\beta
	\nonumber\\
	&\geq \frac{d_j^2 (1-l_{j}/l_{j+1})}{8(2a_j+1)} \int_{K l_{j+1}}^{t-r} \alpha^{-1}
	\log^{2a_{j}+1} \left( \frac{\alpha}{K l_{j+1}} \right) d \alpha\nonumber\\
	&\geq \frac{d_j^2 (1-l_{j}/l_{j+1})}{8(a_{j+1})^2}
	\log^{2a_j+2} \left( \frac{t-r}{K l_{j+1}} \right).
	\label{dai4.71}
\end{align}

Solving (\ref{dai4.69}) yields
\begin{align}
	a_{j}=2^{j+1}-2.\label{dai4.72}
\end{align}
Since $(1-l_{j}/l_{j+1})=2^{-(j+1)}/l_{j+1} \geq 2^{-(j+2)}$, we have
\begin{align}
	\frac{1-l_{j}/l_{j+1}}{8(a_{j+1})^2}
	\geq \frac{ 2^{-(j+2)}}{2^3 \cdot 2^{2j+4}}
	= \frac{1}{2^{3j+9}}.\label{dai4.73}
\end{align}
Therefore, from (\ref{dai4.71}), (\ref{dai4.72}) and (\ref{dai4.73}),
 (\ref{dai4.68}) holds for all natural numbers.
\par
We get from (\ref{dai4.70})
\begin{align}
	\log d_{j+1}=2^{j+1} \log d_0-(\log 2) \sum_{k=0}^{j}
	\left\{ (3(j-k)+9)2^{k} \right\}.
	\nonumber
\end{align}
We obtain
\begin{align}
	d_{j}=\exp\left\{ 
	2^{j} \left(\log d_{0} -(\log 2)
	\sum_{k=0}^{j-1} \frac{(3(j-k)+9)2^{k}}{2^j}
	\right)
	\right\}.\label{dai4.74}
\end{align}
The sum part in (\ref{dai4.74}) converges as $j \to \infty$ by the d'Alembert's criterion.
Hence, there exists a constant $q$ such that it holds
\begin{align}
	d_{j} \geq \exp \left\{ 2^{j} \log \left( E_{0} e^q \e \right) \right\}.\label{dai4.75}
\end{align}
Since $l_{j} \leq 2$, we get
from (\ref{dai4.68}), (\ref{dai4.72}) and (\ref{dai4.75})
\begin{align}
	\sqrt{(t+r)(t-r)}\bar{u}(r,t)\geq
	\exp \left\{ 2^j J(r,t) \right\} \log^{-2} \left( \frac{t-r}{2K} \right)
	\quad \mbox{in} \ \Sigma_{\infty},
	\label{dai4.76}
\end{align}
where
\begin{align}
	&J(r,t):= \log \left(  \e \left\{ 
	B^{-1} \log\left( \frac{t-r}{2K} \right) \right\}^2\right)
	\quad \mbox{and} \quad 
	B:=E_0^{-1/2} e^{-q/2}.
	\label{dai4.77}
\end{align}
\par
We take $\e_{0}>0$ so small that
\begin{align}
	B\e_{0}^{-1/2} \geq \log(2K).\label{dai4.78}
\end{align}
For a fixed $\e \in [0,\e_0)$, we suppose that $T$
satisfies
\begin{align}
	T \ge \exp \left( 4B \e^{-1/2} \right).
	\label{dai4.79}
\end{align}
Next, we take $\tau>0$ so that
\[
T>\tau > \exp \left( 2B \e^{-1/2} \right)\quad (>2K).
\]
From (\ref{dai4.78}) and (\ref{dai4.79}), it follows  that
\begin{align}
	\tau > 2K \exp(B \e^{-1/2}).\label{dai4.80}
\end{align}
We get from (\ref{dai4.77}) and (\ref{dai4.80})
\begin{align}
	J(0, \tau) =\log \left( \e 
	\left\{ B^{-1}\log \frac{\tau}{2 K} \right\}^2 \right)>0.
	\label{dai4.81}
\end{align}
Since $(0, \tau)\in \Sigma_{\infty}$, from (\ref{dai4.76})
and (\ref{dai4.81}), we get $u(0,\tau) \to \infty$ $(j \to \infty)$. 
This completes the proof.
\hfill$\Box$

\subsection{Proof of Theorem \ref{T.2.2}}
\par
We divide the proof of Theorem \ref{T.2.2} into two cases,
$1<p<2$ and $p=2$.
First, we shall handle the sub-critical case.\\
\par\noindent
{\bf Proof of Theorem \ref{T.2.2} with \mbox{\boldmath{$1<p<2$}}.}
Due to the assumption on the initial data in Theorem \ref{T.2.2},
$f(x)\ge0$ $(\not\equiv0)$ and $f(x)+g(x)\equiv0$,
we have
\begin{align*}
	F(0)>0,\quad F'(0)=0.
\end{align*}
For the key inequality, we employ the following proposition.
\begin{prop}{\label{prop4.1}}
Let $1<p<2$.
Suppose that the assumptions 
in Theorem \ref{T.2.2} are fulfilled. Then, there exists a positive
constant $C_{*}=C_{*}(f,g,p,k)$ such that $F(t)=\int_{\R^2} u(x,t) dx$ satisfies
\begin{align}
	F''(t) \geq C_{*}\e^p t^{2-3p/2}
	\quad \mbox{for} \ t \geq k.\label{dai4.22}
\end{align}	
\end{prop}
{\bf Proof.}
It follows from Lemma 2.2 in \cite{YZ06} that
\begin{align}
	F_{1}(t) \geq  \frac{1}{2} \left( 1-e^{-2k}\right)
	\int_{\R^2} \e f(x)  \phi_{1} (x) dx
	\quad \mbox{for} \ t  \geq k,\label{dai4.37}
\end{align}
where
\begin{align*}
	\phi_{1}(x) = \int_{S^{1}} e^{x \cdot \omega} d \omega \quad \mbox{and} \quad 
	F_{1}(t) =e^{-t} \int_{\R^2} u(x,t) \phi_{1}(x) dx.
\end{align*}
From (2.4) and (2.5) in \cite{YZ06}, we obtain
\begin{align}
	F''(t) \geq C(t+k)^{2-3p/2} |F_{1}(t)|^{p} \quad \mbox{for} \ t \geq 0.
	\label{dai4.38}
\end{align}
From (\ref{dai4.37}) and (\ref{dai4.38}), we obtain (\ref{dai4.22}). This completes the proof.
\hfill$\Box$\\
In the sub-critical case,
the following basic lemma is useful.
\begin{lem}[Takamura \cite{Takamura15}]\label{lem4.3}
Assume that (\ref{dai4.11}) is replace by
\begin{align}
	& F(0)>0,\ F'(0)=0,\nonumber
\end{align}
and additionally that there is a $t_0>0$ such that  
\begin{equation}\label{dai4.27}
	F(t_0)\ge 2F(0).
\end{equation}
Then, the conclusion of Lemma \ref{lem4.1} is changed to that 
there exists a positive constant $\widetilde{D_0}=\widetilde{D_0}(p,a,q,B)$ such that
\begin{equation}
	T<2^{2/M}T_2\nonumber
\end{equation}
holds provided
\begin{equation}
	T_2:=\max\left\{T_0,t_0,k\right\}\ge \widetilde{D_0} A^{-(p-1)/(2M)}.
	\label{dai4.29}
\end{equation}
\end{lem}
This is exactly Lemma 2.2 in \cite{Takamura15}, so that we shall omit the proof here.
\par
Integrating (\ref{dai4.22}) in $[k,t]$, we have
\begin{align*}
	F'(t) \geq \frac{C_{*}}{3-3p/2} \e^p 
	(t^{3-3p/2}-k^{3-3p/2})+F'(k)
\end{align*}
for $t \geq k$ because of $1 < p < 2$.
Note that $F'(k)\geq0$ follows from $F''(t) \geq 0$ for $t \geq 0$ and $F'(0) =0$.
Hence we obtain that
\begin{align*}
	F'(t) \geq \frac{C_{*}(1-2^{-3+3p/2})}{3-3p/2} \e^p 
	t^{3-3p/2} \quad \mbox{for} \ t \geq 2k.
\end{align*}
Integrating this inequality in $[2k,t]$ together with $F(0) > 0$, we get
\begin{align}
	F(t) \geq D_{3} \e^p t^{4-3p/2} \quad \mbox{for} \ t\geq 4k,
	\label{dai4.31}
\end{align}
where
\begin{align*}
	D_{3}:=
	\frac{C_{*}(1-2^{-3+3p/2})(1-2^{-4+3p/2})}
	{( 3-3p/2) (4-3p/2)}>0.
\end{align*}
From $2F(0)=2\| f \|_{L^{1}(\R^2)}\e$ and (\ref{dai4.31}), (\ref{dai4.27}) in Lemma \ref{lem4.3} is fulfilled with
\begin{align*}
	t_{0}:=D_{4} \e^{-(p-1)/(4-3p/2)},
\end{align*}
if $t_{0} \geq 4k$, where
\begin{align*}
	D_{4}:=\left\{  2 \| f \|_{L^{1}(\R^2)} D_{3}^{-1} \right\}^{(4-3p/2)^{-1}}.
\end{align*}
We are now in a position to apply our result  here to Lemma 
\ref{lem4.3} with special choices
on all positive constants except for $T_{0}$ as
\begin{align*}
	A=D_{3} \e^p,\quad B=\pi^{1-p},\quad a=4-\frac{3}{2}p,\quad
	q=3(p-1)
\end{align*}
which imply that (\ref{dai4.8}) yields
\begin{align*}
	M=\frac{p-1}{2}a-\frac{q}{2}+1=\frac{\gamma(p,4)}{4}>0.
\end{align*}
If we set
\begin{align*}
	T_{0}:=\widetilde{D_{0}}A^{-(p-1)/(2M)}=\widetilde{D_{0}}D_{3}^{-2(p-1)/\gamma(p,4)} 
	\e^{-2p(p-1)/\gamma(p,4)},
\end{align*}
then we find that there is an $\e_{0}=\e_{0}(f,g,n,p,k)>0$ such that
\begin{align*}
	T_{0} \geq \max \{t_{0}, 4k \} \quad 
	\mbox{for} \ 0<\e \leq \e_{0}
\end{align*}
because of $2p/\gamma(p,4) < 1/(4-3p/2)$.
This means that
$T_{2}=T_{0}$ in (\ref{dai4.29}). Therefore, the conclusion of Lemma \ref{lem4.3} implies 
\begin{align*}
     T < 2^{2/M} T_{2} = D_{5}\e^{-2p(p-1)/\gamma(p,4)}
	\quad \mbox{for} \ 0 <\e \leq \e_{0},
\end{align*}
where
\begin{align*}
	D_{5}:=2^{8/\gamma(p,4)} \widetilde{D_{0}} D_{3}^{-2(p-1)/\gamma(p,4)}>0.
\end{align*}
This completes the proof.
\hfill$\Box$\\
\par\noindent
{\bf Proof of Theorem \ref{T.2.2} with \mbox{\boldmath{$p=2$}}.}
Since 
\[
f(x)+g(x)\equiv0\quad\mbox{and}\quad\int_{\R^2}f(x)dx<0,
\]
by Lemma \ref{lem4.5}, there exist positive constants $\widetilde{E_0}$ and $\widetilde{K} \geq 1$
such that 
\begin{align}
	u_{L}(x,t) \geq \frac{\widetilde{E_0}  }{(t+r)^{1/2} (t-r)^{3/2}} \quad
	\mbox{for} \ t-r \geq \widetilde{K}.\nonumber
\end{align}
For $t-r \geq \widetilde{K}$, we get from (\ref{dai4.61}) that
\begin{align}
	\bar{u}(r,t) \geq \overline{u^{0}}(r,t) \geq 
	\frac{\widetilde{E_0} \e }{(t+r)^{1/2}(t-r)^{3/2}}.\label{dai4.90}
\end{align}
We define the following domains:
\begin{align}
	&\widetilde{\Sigma}_{j}:=\left\{
	(r,t) \in [0,\infty)^2 \ | \ t-r \geq 3 \widetilde{K} l_{j}
	\right\} \quad (j=0,1,2,\cdots),\nonumber\\
	&\widetilde{\Sigma}_{\infty}:=\left\{
	(r,t) \in [0,\infty)^2 \ | \ t-r \geq 6 \widetilde{K}
	\right\}.\nonumber
\end{align}
In the same way as to obtain (\ref{dai4.67}), for $t-r  \geq \widetilde{K}$, we get
\begin{align}
	\bar{u}(r,t) \geq \frac{1}{\sqrt{(t+r)(t-r)}}
	\int\int_{\widetilde{Q_0}(r,t)} (1+\tau)^{-1} \lambda |\bar{u}|^2 d\lambda d\tau,
	\label{dai4.91}
\end{align}
where
\begin{align}
	\widetilde{Q_0}(r,t):=\left\{
	(\lambda,\tau) \in [0,\infty)^2 \ | \  \widetilde{K} \leq \tau - \lambda,\
	0 \leq \lambda \leq t-r-\tau
	\right\}.\nonumber
\end{align}
For $(r,t) \in \widetilde{\Sigma_0}$, we set
\begin{align}
	S(r,t):=\left\{
	(\lambda,\tau) \in [0,\infty)^2 \ \Big| \ \frac{5}{2}\widetilde{K} \leq \tau + \lambda \leq
	t-r, \ \widetilde{K} \leq \tau -\lambda \leq \frac{5}{4} \widetilde{K}
	\right\}. \nonumber
\end{align}
For $(r,t) \in \widetilde{\Sigma_0}$, we have $S(r,t) \subset \widetilde{Q_0}(r,t)$.
Substituting (\ref{dai4.90}) into (\ref{dai4.91}) and changing the variables by (\ref{dai2.3}), we get
\begin{align}
	&\sqrt{(t+r)(t-r)} \bar{u}(r,t)\nonumber\\
	&\geq \int \int_{S(r,t)} (1+\tau)^{-1} \lambda |\bar{u}|^2 d \lambda d \tau
	\nonumber\\
	&\geq \frac{1}{2} \int_{5\widetilde{K}/2}^{t-r} d \alpha
	\int_{\widetilde{K}}^{5\widetilde{K}/4} 
	\left(1+ \frac{\alpha+\beta}{2} \right)^{-1} 
	\left( \frac{\alpha-\beta}{2} \right) 
	\left(\frac{\widetilde{E_0} \e}{\alpha^{1/2} \beta^{3/2}} \right)^2
	d \beta.\label{dai4.92}
\end{align}
Since $\alpha+\beta \leq 2 \alpha$, $\displaystyle \alpha- \beta \geq \alpha-5k/4 \geq \alpha/2 $ and $\widetilde{K} \geq 1$, we get from (\ref{dai4.92})
\begin{align}	
	\sqrt{(t+r)(t-r)} \bar{u}(r,t)
	\geq& \frac{(\widetilde{E_0} \e )^2}{2^4}
	\int_{5\widetilde{K}/2}^{t-r} \alpha^{-1} d \alpha 
	\int_{\widetilde{K}}^{5\widetilde{K}/4} \beta^{-3/2} d\beta\nonumber\\
	\geq & 
	E_{0}^{*}\e^2 \log \left( \frac{t-r}{3 \widetilde{K}} \right)
	\quad \mbox{in} \ \widetilde{\Sigma_0},\label{dai4.93}
\end{align}
where $\displaystyle E_{0}^{*}:=\widetilde{E_0}^2/\left(2^9 \widetilde{K}^{1/2} \right)$.
Analogously to the proof of Theorem \ref{T.2.1} with $p=2$, we obtain
from (\ref{dai4.93})
\begin{align}
	\bar{u}(r,t) \geq \frac{\widetilde{d_j}}{\sqrt{(t+r)(t-r)}}
	\log ^{\widetilde{a_j}} \left( \frac{t-r}{3 \widetilde{K} l_{j}}\right)
	\quad \mbox{in} \ \widetilde{\Sigma_j},\nonumber
\end{align}
where
\begin{align}
	&\widetilde{a_0}=1,\quad \widetilde{a_{j+1}}=2\widetilde{a_{j}}+2,\nonumber\\
	&\widetilde{d_0}=E_0^{*} \e^2,\quad \widetilde{d_{j+1}}
	=\frac{\widetilde{d_j}^2}{3 \cdot 2^{3j+9}}.\nonumber
\end{align}
This is the same form as (\ref{dai4.68}), (\ref{dai4.69}) and (\ref{dai4.70}). Hence, we see that there exists a constant $\widetilde{q}$ such that
\begin{align*}
	\widetilde{d_j} \geq \exp \left\{ 2^{j} \log 
	\left( E_{0}^{*} e^{\widetilde{q}} \e^2 \right) \right\}.
\end{align*}
Since $l_j \leq 2$, we get
\begin{align*}
	\sqrt{(t+r)(t-r)} \bar{u}(r,t)
	\geq \exp \left\{ 2^{j} J(r,t) \right\} \log^{-2} 
	\left( \frac{t-r}{6 \widetilde{K}} \right) \quad \mbox{in} \ \widetilde{\Sigma_\infty},
\end{align*}
where
\begin{align*}
	J(r,t) := \log \left\{ \e^2 \left\{ \widetilde{B}^{-1}
	\log \left( \frac{t-r}{6 \widetilde{K}} \right)  \right\}^3\right\}
	\quad \mbox{and} \quad \widetilde{B}:=(E_{0}^{*})^{-1/3}e^{-\widetilde{q}/3}.
\end{align*}
In the same way as in the proof of Theorem \ref{T.2.1} with $p=2$,
we get the desired estimates. 
The proof is now completed.
\hfill$\Box$

\section*{Acknowledgement}
\par
This work started when the first author was a master course student in Future University Hakodate,
the third author was working in Future University Hakodate,
and the fourth author was working in Muroran Institute of Technology.
The third author has been partially supported by
Special Research Expenses in FY2017, General Topics (No.B21), Future University Hakodate,
also by the Grant-in-Aid for Scientific Research (B) (No.18H01132) and (C) (No.15K04964), 
Japan Society for the Promotion of Science.
Finally all the authors are grateful to the referee for precise reading and useful comments.


\bibliographystyle{plain}

\end{document}